\documentclass[11pt,article]{article}

\usepackage{amsmath}
\usepackage{amsthm}
  \usepackage{paralist}
  \usepackage{graphics} 
  \usepackage{epsfig} 
\usepackage{graphicx}
\usepackage{caption}
\usepackage{subcaption}
\usepackage{epstopdf}
 \usepackage[colorlinks=true]{hyperref}
 \usepackage{multirow}
\input{amssym.tex}

\newtheorem{theorem}{Theorem}[section]
\newtheorem{corollary}{Corollary}
\newtheorem{lemma}[theorem]{Lemma}
\newtheorem{proposition}{Proposition}

 \numberwithin{equation}{section}
\newtheorem{remark}{Remark}

\newcommand{\keywords}

\def\bc{\begin{center}}       \def\ec{\end{center}}
\def\ba{\begin{array}}        \def\ea{\end{array}}
\def\be{\begin{equation}}     \def\ee{\end{equation}}
\def\bea{\begin{eqnarray}}    \def\eea{\end{eqnarray}}
\def\beaa{\begin{eqnarray*}}  \def\eeaa{\end{eqnarray*}}

\def\mathbb{\Bbb}

\begin{document}

\title{\bf Global existence and steady states of a two competing species Keller--Segel chemotaxis model}
\author{Jia Hu \thanks{ {\tt jiahu@2011.swufe.edu.cn}.}, Qi Wang \thanks{ {\tt qwang@swufe.edu.cn}.}, Jingyue Yang   \thanks{ {\tt yjy@2011.swufe.edu.cn}.}, Lu Zhang \thanks{ {\tt lz@2011.swufe.edu.cn}.}\\
Department of Mathematics\\
Southwestern University of Finance and Economics\\
555 Liutai Ave, Wenjiang, Chengdu, Sichuan 611130, China
}

\date{}
\maketitle

\abstract
We study a one--dimensional quasilinear system proposed by J. Tello and M. Winkler \cite{TW} which models the population dynamics of two competing species attracted by the same chemical.  The kinetic terms of the interacting species are chosen to be of the Lotka--Volterra type and the boundary conditions are of homogeneous Neumann type which represent an enclosed domain.  We prove the existence of globally bounded classical solutions for the fully parabolic system.  Moreover, we establish the existence of nonconstant positive steady states through bifurcation theory.  The stability or instability of the bifurcating solutions is also investigated rigorously.  Our results indicate that small intervals support stable monotone positive steady states and large intervals support nonmonotone steady states.  Finally, we perform extensive numerical studies to verify our theoretical results.  Our numerical simulations demonstrate the formation of stable positive steady states and time--periodic solutions with interesting spatial structures.


\section{Introduction}\label{section1}

This paper is devoted to study the global solutions and nonconstant positive steady states of the following reaction--advection--diffusion system,
\begin{equation}\label{11}
\left\{
\begin{array}{ll}
u_t=(d_1u_x-\chi u w_x)_x+\mu_1(1-u-a_1v)u,&x \in (0,L),t>0,      \\
v_t=(d_2v_x-\xi v w_x)_x+\mu_2(1-a_2u-v)v,&x \in (0,L),t>0,      \\
w_t=w_{xx}-\lambda w +u+v,&x \in (0,L),t>0,      \\
u_x(x,t)=v_x(x,t)=w_x(x,t)=0,&x=0,L,t>0,\\
u(x,0)=u_0(x),v(x,0)=v_0(x),w(x,0)=w_0(x),&x\in(0,L),
\end{array}
\right.
\end{equation}
where $u$, $v$ and $w$ are functions of space $x$ and time $t$.  $d_i$, $\mu_i$, $a_i$, $i=1,2$, and $\lambda$ are positive constants; $\chi$ and $\xi$ are assumed to be nonnegative constants. (\ref{11}) is the one--dimensional version of the following model with $\tau=1$
\begin{equation}\label{12}
\left\{
\begin{array}{ll}
u_t=\nabla \cdot (d_1 \nabla u-\chi u \nabla w)+\mu_1(1-u-a_1v)u,&x \in \Omega,t>0, \\
v_t=\nabla \cdot(d_2\nabla v-\xi v \nabla w)+\mu_2(1-a_2u-v)v,&x \in \Omega,t>0, \\
\tau w_t=\Delta w-\lambda w +u+v,&x \in \Omega,t>0,\\
\frac{\partial u}{\partial \textbf{n}}=\frac{\partial v}{\partial \textbf{n}}=\frac{\partial w}{\partial \textbf{n}}=0,&x\in\partial \Omega,t>0,\\
u(x,0)=u_0(x),v(x,0)=v_0(x),w(x,0)=w_0(x), &x\in \Omega,
\end{array}
\right.
\end{equation}
which was considered by J. Tello and M. Winkler \cite{TW} to study the spatial--temporal behaviors of two competing species attracted by the same chemical.  It is assumed in \cite{TW} that $\Omega \subset \mathbb{R}^N$, $N\geq1$, is a bounded domain with smooth boundary $\partial \Omega$.  $u(x,t)$ and $v(x,t)$ represent the population densities of two competing species at space--time location $(x,t)\in \Omega\times \mathbb{R}^+$, while $w(x,t)$ denotes the concentration of the attracting chemical.  It is assumed that both species $u$ and $v$ direct their movements chemotactically along the gradient of the chemical concentration over the habitat.  This is modeled by taking both $\chi>0$ and $\xi>0$.  Biologically, $\chi$ and $\xi$ measure the strength of chemical attraction to species $u$ and $v$ respectively.  The kinetics of the species are assumed to be of the classical Lotka--Volterra type, where $\mu_1$, $\mu_2$ measure the intrinsic growth rates of $u$, $v$ and $a_1$, $a_2$ interpret the strength of interspecific competition.  Moreover, the chemical is produced by both species at the same rate with no saturation effect and is consumed by certain enzyme in the environment at rate $\lambda$ meanwhile.

System (\ref{12}) has a positive constant steady state
\begin{equation}\label{13}
(\bar u,\bar v,\bar w)=\Big( \frac{1-a_1}{1-a_1a_2}, \frac{1-a_2}{1-a_1a_2},\frac{2-a_1-a_2}{\lambda (1-a_1a_2)} \Big),
\end{equation}
provided that
\begin{equation}\label{14}
0\leq a_1<1,0\leq a_2<1.
\end{equation}

J. Tello and M. Winkler investigated the global existence of parabolic--parabolic--elliptic system (\ref{12}) and the global asymptotic stability of the constant equilibrium $(\bar u,\bar v,\bar w)$.  Their results can be summarized as follows:
\begin{theorem}\label{theorem11}(Tello, Winkler \cite{TW})
Let $\Omega\subset \mathbb{R}^N$, $N\geq1$, be a bounded domain with smooth boundary $\partial \Omega$.  Assume that $\tau=0$, $d_1=d_2=1$ and condition (\ref{14}) is satisfied; moreover, suppose that the positive constants $\chi$, $\xi$, $\mu_1$ and $\mu_2$ satisfy the following condition
\begin{equation}\label{15}
2(\chi+\xi)+a_1\mu_2<\mu_1,~2(\chi+\xi)+a_2\mu_1<\mu_2,
\end{equation}
then for any positive initial data $(u_0,v_0,w_0)\in C^0(\bar \Omega)\times C^0(\bar \Omega)\times C^0(\bar \Omega)$, the solution $(u(x,t),v(x,t),w(x,t))$ of (\ref{12}) is global and bounded and $(\bar u,\bar v,\bar w)$ given by (\ref{13}) is a global attractor of the parabolic--parabolic--elliptic system of (\ref{12}) in the sense
\[\Vert u(\cdot,t)-\bar u \Vert_{L^\infty}+\Vert v(\cdot,t)-\bar v \Vert_{L^\infty} \rightarrow 0,t\rightarrow \infty.\]
\end{theorem}
Theorem \ref{theorem11} suggests that for any given $\mu_1,\mu_2>0$, small chemo--attraction rates $\chi$, $\xi$ and competition rates $a_1$, $a_2$ stabilize the positive equilibrium $(\bar u,\bar v,\bar w)$, which is globally asymptotically stable.  Therefore the stationary system of (\ref{12}) has no positive solutions other than $(\bar u,\bar v,\bar w)$ under condition (\ref{15}).  From the viewpoint of mathematical modeling, it is interesting to investigate nonconstant positive steady states of the two species chemotaxis model (\ref{12}).  Moreover, it is important to study the evolution of species population distributions through interspecific competition and chemotaxis processes.  It is the motivation of this paper to study the range of parameters $\chi$ and $\xi$ in which the one--dimensional model (\ref{11}) admits nonconstant positive solutions.  In particular, we study the existence and stability of nonconstant positive steady states of system (\ref{11}).

There are many experimental examples on the effects of chemotaxis on the dynamics of two microbial populations competing for single or multiple rate-limiting nutrient(s) in a confined system, for example by \cite{AT, TMK, VL} etc.  D. Lauffenburger \emph{et al}. \cite{KDL, L, LAK, LC} initiated theoretical analysis of the effects of cell motility and chemotaxis on the population dynamics of single or two competing microbial populations with confined growth in a tubular reactor, which is supplied with a single diffusible growth--limiting nutrient entering at one end of the tube.  Their results suggest that in such nonmixed systems, cell motility and chemotaxis properties are the determining factors in governing population dynamics.

Some other two--species chemotaxis systems closely related to (\ref{12}) are proposed and studied by various authors.   For $\mu_1=\mu_2=\tau=0$ and $\lambda w$ being replaced by a positive constant, E. Espejo \emph{et al}. \cite{ESV} investigated  simultaneous finite--time blow--ups of (\ref{12}) when $\Omega$ is a circle in $\mathbb{R}^2$.  For $\mu_1=\mu_2=\lambda=\tau=0$, P. Biler \emph{et al}. \cite{BEG} obtained the blow--up properties of (\ref{12}) with $\Omega=\mathbb{R}^N$, $N\geq2$.  Similar blow--up mechanisms, in particular related to the initial data size, have been studied in \cite{CEV} and \cite{CEV2} for $\Omega=\mathbb{R}^2$.  X. Wang and Y. Wu \cite{WW} studied the qualitative behaviors of two--competing species chemotaxis models with cellular growth kinetics different from those in (\ref{11}).  In particular, they investigated the large time behaviors of global solutions and the existence of nonconstant positive steady states.  Moreover, the effect of cell motility and chemotaxis on population growth has also been analyzed.  D. Horstmann \cite{Ho} proposed some multi--species chemotaxis models with attraction and repulsion between interacting species.  A linearized stability analysis of the constant steady state is performed and the existence of Lyapunov functional is also obtained for some of the proposed new models.  P. Liu \emph{et al}. \cite{LSW} investigated the pattern formation of an attraction--repulsion Keller--Segel system with two chemicals and one species.  We also want to mention that the two--species Lotka--Volterra competition system with no chemical is recently investigated by Q. Wang \emph{et al}. in \cite{WGY}.

It seems necessary to point out that none of the works above on chemotaxis except \cite{TW} considered Lotka--Volterra competition dynamics, which are usually applied to describe population dynamics in ecological systems.  Cellular chemotaxis and population dynamics are on different spatial and temporal scales, therefore it lacks the rationale to act Lotka--Volterra on the chemotactic microbial populations.  However, our results and analysis on model (\ref{11}) carry over to $3\times 3$ system with different kinetics modeling cellular competitions.  It is our focus to study the dynamics and pattern formation of (\ref{11}) due to the effect of chemotaxis.  J. Tello and M. Winkler considered the smallness condition on $\chi$ and $\xi$ when $(\bar u,\bar v,\bar w)$ is a global attractor of (\ref{12}) with $\tau=0$.  Our work complements theirs by studying its nonconstant positive solutions and steady states over one--dimensional finite intervals, in particular those with interesting patterns.

In this paper,  we study positive solutions to system (\ref{11}) and its stationary system.  We are concerned with the global existence of system (\ref{11}) and the formation of its spatially inhomogeneous steady states when the smallness condition (\ref{15}) on $\chi$ and $\xi$ is relaxed.  Our results state that (\ref{11}) admits a unique pair of classical solution $(u,v,w)$ which is uniformly bounded in $L^\infty$ for all $t\in(0,\infty)$--see Theorem \ref{theorem25}.  In Section \ref{section3}, we establish the existence and stability of nonconstant positive steady states by rigorous bifurcation analysis--see Theorem \ref{theorem31} and Theorem \ref{theorem32}.  In loose terms, large chemotaxis rate $\chi$ (or large $\xi$) destabilizes the homogeneous positive steady state and nonconstant positive steady states can emerge and form stable patterns.  In Section \ref{section4}, we present some numerical simulations to illustrate the formation of self--organized spatial patterns of (\ref{11}) that have complicated and interesting structures.  Finally, we include our concluding remarks and propose some interesting questions for future studies in Section \ref{section5}.

\section{Existence of global--in--time solutions}\label{section2}
In this section, we study the global existence of uniformly bounded and classical solutions $(u,v,w)$ to (\ref{11}).  We shall first apply the well--known results of Amann \cite{A,A1} to obtain the local existence, then we obtain the global existence by establishing the $L^\infty$--bounds of $(u,v,w)$.

\subsection{Local existence}
We present the existence and uniqueness of the positive local classical solutions of (\ref{11}) as follows.
\begin{theorem}\label{theorem21}
Assume that the constants $a_i$, $d_i$, $\mu_i>0$ for $i=1,2$.  Suppose that the initial data $u_0$, $v_0$ and $w_0\in H^1(0,L)$ satisfy $u_0$, $v_0$ and $w_0\geq0$, $\not \equiv0$ on $[0,L]$.  Then for any $\chi$, $\xi\in \mathbb{R}$, the following statements hold true:

(i)  There exists a unique solution $(u(x,t),v(x,t),w(x,t))$ of (\ref{11}) which is nonnegative on $[0,L]\times [0,T_{\max})$ with $0<T_{\max} \leq \infty$ such that $(u(\cdot,t),v(\cdot,t),w(\cdot,t)) \in C([0,T_{\max}),H^1(0,L)\times H^1(0,L)\times H^1(0,L))$ and $(u,v,w)\in C^{2+\alpha,2+\alpha,1+\alpha}_{loc}([0,L]\times (0,T_{\max}))$ for any $0<\alpha<\frac{1}{4}$.

(ii)  If $\sup_{s\in(0,t)}\Vert (u,v,w)(\cdot,s) \Vert_{L^\infty}$ is bounded for $t\in(0,T_{\max})$, then $T_{\max}=\infty$, i.e., $(u,v,w)$ is a global solution to (\ref{11}).  Furthermore, $(u,v,w)$ is a classical solution and $(u,v,w) \in C^\alpha((0,\infty),C^{2(1-\beta)}([0,L])\times C^{2(1-\beta)}([0,L])\times C^{2(1-\beta)}([0,L]))$ for any $0\leq \alpha \leq \beta \leq 1$.
\end{theorem}

\begin{proof}
We write (\ref{11}) into the abstract form
\begin{equation}\label{21}
\left(
\begin{array}{c}u\\
v\\
w
\end{array}
\right)_t=
\left[ \mathcal{A}_0 \left(
\begin{array}{c}u\\
v\\
w
\end{array}
\right)_x \right]_x
+\left(
\begin{array}{c}
\mu_1(1-u-a_1v)u\\
\mu_2(1-a_2u-v)v\\
-\lambda w+u+v
\end{array}
\right),
\end{equation}
where
\begin{equation*}\mathcal{A}_0 =\begin{pmatrix}
 d_1  &  0& -\chi u  \\
  0    & d_2& -\xi v\\
  0       & 0    &  1
  \end{pmatrix}.
  \end{equation*}
System (\ref{21}) is normally parabolic since all the eigenvalues of $\mathcal{A}_0$ are positive.  Then \emph{(i)} follows from Theorem 7.3 and Theorem 9.3 of \cite{A}.  Moreover, \emph{(ii)} follows from Theorem 5.2 in \cite{A1} since (\ref{21}) is a triangular system.  The nonnegativity of $(u,v,w)$ follows from parabolic Maximum Principles--see \cite{LSU}.
\end{proof}

\subsection{A--prior estimates and global existence}
We now proceed to establish the $L^\infty$--bounds of $(u,v,w)$ under the same conditions in Theorem \ref{theorem21}.  The global existence of (\ref{11}) is a consequence of several lemmas.
\begin{lemma}\label{lemma22}
Let $\mu_1$, $\mu_2>0$ and $(u,v,w)$ be the unique nonnegative solution of (\ref{11}), then there exists a constant $C>0$ dependent on $\Vert (u_0,v_0,w_0 )\Vert_{L^1}$ and $L$ such that
\begin{equation}\label{22}
\Vert u(\cdot,t) \Vert_{L^1(0,L)}+\Vert v(\cdot,t) \Vert_{L^1(0,L)}+\Vert w(\cdot,t) \Vert_{L^1(0,L)}\leq C,~\text{for all}~t\in(0,T_{\max}).
\end{equation}
\end{lemma}

\begin{proof}
We integrate the $u$--equation of (\ref{11}) over $(0,L)$ and have that
\[\frac{d}{dt} \int_0^L u(x,t)dx=\int_0^L \mu_1(1-u-a_1v)udx\leq\int_0^L \mu_1(1-u)udx,\]
then it follows from the Gronwall's lemma that
\[\int_0^L u(x,t)dx\leq e^{-\mu_1 t}\int_0^L u_0(x) dx+L;\]
similarly we can show
\[\int_0^L v(x,t)dx\leq e^{-\mu_2 t}\int_0^L v_0(x) dx+L.\]
Moreover, integrating the $w$--equation over $(0,L)$, we easily see that $\Vert w(\cdot,t) \Vert_{L^1(0,L)}$ is uniformly bounded for all $t\in(0,\infty)$.  This completes the proof of Lemma \ref{lemma22}.
\end{proof}

Lemma \ref{lemma22} gives the $L^1$--bounds of $u,v$ and $w$.  To obtain their $L^\infty$--bounds, we shall see that it is necessary to obtain the boundedness of $\Vert w_x(\cdot,t) \Vert_{L^\infty}$.  For this purpose, we first convert the $w$--equation into the following abstract form
\begin{equation}\label{23}
\begin{split}
w(\cdot,t)=&\left.e^{ (\Delta-1)t}w_0+\int_0^t e^{ (\Delta-1)t} \big((1-\lambda) w(\cdot,s)+u(\cdot,s)+v(\cdot,s) \big)ds, \right.
\end{split}
\end{equation}
where $\Delta=\frac{d^2}{dx^2}$.  To estimate $w(x,t)$ in (\ref{23}), we apply the well--known smoothing properties of operator $-\Delta+1$ and estimates between the linear analytic semigroups generated by $\{e^{t\Delta}\}_{t\geq 0}$, for which \cite{H, HW, Winkler} are good references.  For example, we have from Lemma 1.3 with $N=1$ in \cite{Winkler} that for all $1\leq p \leq q \leq \infty$, there exists a positive constant $C$ dependent on $\mu_1$, $\mu_2$, and $\Vert w_0 \Vert _{W^{1,q}(0,L)}$ such that
\begin{equation}\label{24}
 \Vert w(\cdot,t) \Vert _{W^{1,q}} \!\leq \!C\!\!\left(1\! + \!\int_0^t\!\!e^{-\nu(t-s)} (t-s)^{-\frac{1}{2}-\frac{1}{2}(\frac{1}{p}-\frac{1}{q})}   \Vert w(\cdot,s)\!+\!u(\cdot,s)\!+\!v(\cdot,s) \Vert_{L^p}   ds\!\right),
\end{equation}
for any $t\in(0,T)$, all $T\in(0,\infty]$, where $\nu$ is the first Neumann eigenvalue of $-\Delta$.  We have the following a--prior estimate on $\Vert w(\cdot,t)\Vert_{W^{1,q}(0,L)}$.
\begin{lemma}\label{lemma23}
Assume the same conditions on $(u_0,v_0,w_0)$ as in Lemma \ref{lemma22}.  For any $q\in(1,\infty)$, there exists a positive constant $C(q)$ such that
\begin{equation}\label{25}
\Vert w (\cdot,t) \Vert_{W^{1,q}(0,L)} \leq C(q), \forall t\in(0,T_{\max}).
\end{equation}
\end{lemma}

\begin{proof}
Choosing $p=1$ in (\ref{24}), we have that
\begin{equation}\label{26}
\begin{split}
\Vert w (\cdot,t) \Vert_{W^{1,q}(0,L)} \leq C\Big(1+ \int_0^t &e^{-\nu(t-s)} (t-s)^{-\frac{1}{2}-\frac{1}{2}(1-\frac{1}{q})}\cdot\big(\Vert u(\cdot,s) \Vert_{L^1(0,L)}
\\&\left. +\Vert v(\cdot,s) \Vert_{L^1(0,L)}+\Vert w(\cdot,s) \Vert_{L^1(0,L)}\big) ds\right.\Big),
\end{split}
\end{equation}
where $C$ depends on $\Vert (u_0,v_0,w_0) \Vert_{W^{1,q}(0,L)}$.  On the other hand, there exists a constant $C_0>0$ such that for any $q\in(1,\infty)$,
\[\sup_{t\in(0,\infty)}\int_0^te^{-\nu(t-s)} (t-s)^{-1+\frac{1}{2q}} ds<C_0,\]
then we conclude from (\ref{26}) that
\begin{equation}\label{27}
\Vert w (\cdot,t) \Vert_{W^{1,q}(0,L)}\!\! \leq C\Big(\!1+\!\! \sup_{s\in(0,t)}\!\! (\Vert u(\cdot,s) \Vert_{L^1(0,L)}+\Vert v(\cdot,s) \Vert_{L^1(0,L)}+\Vert w(\cdot,s) \Vert_{L^1(0,L)}\!)\Big).
\end{equation}
Finally, (\ref{25}) is an immediate consequence of (\ref{22}) and (\ref{27}).
\end{proof}
Lemma \ref{lemma23} is essential in our proof of the global existence of (\ref{11}).  We want to point out that estimate (\ref{25}) is heavily involved with the space dimension $N=1$ and $q$ can not be arbitrarily large if $N\geq2$.  In higher space--dimensions, the Laplacian may not be sufficient to prevent finite or infinite--time blow--up of (\ref{12}) caused by the chemotaxis.  See our discussions in Section \ref{section5}.
\begin{lemma}\label{lemma24}
For each $p\in(2,\infty)$, there exists a constant $C(p)>0$ such that
\begin{equation}\label{28}
\Vert u(\cdot,t) \Vert_{L^p}+\Vert v(\cdot,t) \Vert_{L^p}\leq C(p),\forall t\in(0,T_{\max}).
\end{equation}
\end{lemma}
\begin{proof}
We shall only show that $\Vert u(\cdot,t) \Vert_{L^p}\leq C(p)$ since $\Vert v(\cdot,t) \Vert_{L^p}\leq C(p)$ can be proved by the same arguments.  For $p>2$, we multiply the first equation of (\ref{11}) by $u^{p-1}$ and integrate it over $(0,L)$ by parts, then it follows from simple calculations that
\begin{eqnarray}\label{29}
\frac{1}{p}\frac{d}{dt}\int_0^L\!\! u^p
\!\!\!\!\!&=&\!\!\!\!\! \int_0^L u^{p-1}u_t =\!\!\int_0^L u^{p-1} (d_1 u_x-\chi u w_x)_x +\!\!\int_0^L\!\! \mu_1(1-u-a_1v)u^p       \nonumber \\
\!\!\!\!\!&\leq&\!\!\!\!\! -\frac{4d_1(p-1)}{p^2}\int_0^L \!\!\vert  (u^\frac{p}{2})_x \vert^2 +\chi\int_0^L\!\!(u^{p-1})_xuw_x -\!\!\frac{\mu_1}{2}\!\!\int_0^L\!\!\!\! u^{p+1}\!+\!C_1   \nonumber \\
\!\!\!\!\!&=&\!\!\!\!\! -\frac{4d_1(\!p\!-\!1\!)}{p^2}\!\!\!\int_0^L\!\!\!\! \vert (u^\frac{p}{2})_x \vert^2 +\!\frac{2(\!p\!-\!1\!)\chi}{p}\!\!\!\int_0^L\!\!\!\!u^\frac{\!p}{2} (u^\frac{p}{2})_x w_x \!-\!\frac{\mu_1}{2}\!\!\int_0^L\!\!\!\! u^{p+1}\!+\!C_1\\ \nonumber
\end{eqnarray}
where $C_1$ is a positive constant that depends on $p$.  On the other hand, it follows from H$\ddot{\text{o}}$lder's and Young's inequality that
\begin{eqnarray}\label{210}
\int_0^L u^\frac{p}{2} (u^\frac{p}{2})_x w_x dx
&\leq& \Vert u^\frac{p}{2} \Vert_{L^{\frac{2(p+1)}{p}}}  \Vert (u^\frac{p}{2})_x \Vert_{L^2}   \Vert w_x\Vert_{L^{2(p+1)}}\nonumber\\
&=& \Vert u \Vert^\frac{p}{2}_{L^{p+1}} \Vert (u^\frac{p}{2})_x \Vert_{L^2} \Vert w_x\Vert_{L^{2(p+1)}} \nonumber\\
&\leq& \frac{2d_1}{p\chi}\Vert (u^\frac{p}{2})_x \Vert^2_{L^2}  +C_2 \Vert u\Vert^p_{L^{p+1}}, \\ \nonumber
\end{eqnarray}
 where we have applied the fact that $\Vert w_x\Vert_{L^{2(p+1)}}\leq C$ due to (\ref{25}).  In light of (\ref{210}), we obtain from (\ref{29}) that
\begin{equation}\label{211}
\frac{1}{p}\frac{d}{dt}\int_0^L u^p dx \leq C_3 \Vert u \Vert^p_{L^{p+1}}-\frac{\mu_1}{2} \Vert u\Vert^{p+1}_{L^{p+1}}+C_1.
\end{equation}
Denoting $y_p(t)=\int_0^L u^p(x,t) dx$, one can apply H$\ddot{\text{o}}$lder's on (\ref{211}) to obtain that
\[y_p'(t)\leq -C_4 y^\frac{p+1}{p}_p(t)+C_5,~y_p(0)=\Vert u_0 \Vert^p_{L^p}.\]
Solving this differential inequality, we conclude that $y_p(t)\leq C(p)$ for all $t\in (0,\infty)$.  Similarly we can show that $\int_0^L v^p(x,t) dx\leq C(p)$.  This finishes the proof of Lemma \ref{lemma24}.
\end{proof}

By taking $p$ large ($\geq 2$) but fixed and $q=\infty$ in (\ref{24}), we can easily obtain the uniform boundedness of $\Vert w_x(\cdot,t) \Vert_{L^\infty}$ thanks to Lemma \ref{lemma24}.
\begin{corollary}\label{corollary1}
Under the conditions in Lemma \ref{lemma22}, there exists a positive constant $C$ such that
\begin{equation}\label{212}
\Vert w(\cdot,t) \Vert_{W^{1,\infty}} \leq C,\forall t\in(0,T_{\max}).
\end{equation}
\end{corollary}

We are now ready to present the following results on the global existence of uniformly bounded and classical positive solutions to (\ref{11}).
\begin{theorem}\label{theorem25}
Let $a_i$, $\mu_i$, $i=1,2$, and $\lambda$ be positive constants.  Then given positive initial data $(u_0,v_0,w_0)\in H^1(0,L)\times H^1(0,L)\times H^1(0,L)$ and any constants $\chi,\xi\in\mathbb{R}$, (\ref{11}) has a unique bounded positive solution $(u(x,t),v(x,t),w(x,t))$ defined on $[0,L]\times [0,\infty)$ such that $(u(\cdot,t),v(\cdot,t),w(\cdot,t)) \in C([0,\infty),H^1(0,L)\times H^1(0,L)\times H^1(0,L))$ and $(u,v)\in C^{2+\alpha,2+\alpha,1+\alpha}_{loc}([0,L]\times [0,\infty))$ for some $\alpha \in (0,\frac{1}{4})$.
\end{theorem}

\begin{proof}According to Part \emph{(ii)} of Theorem \ref{theorem21} and Corollary \ref{corollary1}, we only need to show that $\Vert (u(\cdot,t),v(\cdot,t))\Vert_{L^\infty}$ is uniformly bounded for all $t\in(0,T_{\max})$, then we must have that $T_{\max}=\infty$ and the existence part of Theorem \ref{theorem25} follows.  Moreover, one can apply parabolic boundary $L^p$ estimates and Schauder estimates to show that $u_t,v_t,w_t$ and all spatial partial derivatives of $u$, $v$ and $w$ up to order two are bounded on $[0,L] \times [0,\infty)$, therefore $(u,v,w)$ have the regularities stated in Theorem \ref{theorem25}.

Without loss of our generality, we assume that $\Vert w_x(.,t) \Vert_{L^\infty} \leq 1$ in (\ref{212}).  Otherwise we denote $\chi$ by $\chi C$ without loss of generality in the remaining estimates of the proof, with $C$ being given in (\ref{212}).  Through the same calculations that lead to Lemma \ref{lemma24} and using the fact $u\geq0,v\geq0$, we obtain
\begin{eqnarray}\label{213}
\frac{1}{p} \frac{d}{dt}\!\int_0^L\!\!\!\! u^p dx
\!\!\!\!\!&\leq&\!\!\!\!\!-\frac{4d_1(p-1)}{p^2}\!\int_0^L\!\!\vert (u^\frac{p}{2})_x\vert^2dx+\frac{2(p-1)\chi}{p}\!\int_0^L\!\!u^\frac{p}{2} \vert(u^\frac{p}{2})_x\vert  w_x dx+\mu_1\!\!\int_0^L\!\! u^p dx \nonumber \\
\!\!\!\!\!&\leq&\!\!\!\!\!-\frac{4d_1(p-1\!)}{p^2}\!\int_0^L\! \vert (u^\frac{p}{2})_x \vert^2 dx+\frac{2(p-1)\chi  }{p}\int_0^L u^\frac{p}{2} \vert (u^\frac{p}{2})_x \vert dx+\!\!\mu_1\!\!\!\int_0^L\!\!\!\! u^p dx \nonumber\\
\!\!\!\!\!&\leq&\!\!\!\!\!\!-\frac{4d_1(p-1\!)}{p^2}\!\!\int_0^L\!\!\!\vert(u^\frac{p}{2})_x\! \vert^2 dx\!\!+\!\!\frac{(p-\!1)\chi}{p}\!\!\int_0^L \!\!\!\Big(\frac{p\chi}{2d_1}u^p\!\!+\!\!\frac{2d_1}{p\chi}\vert (u^\frac{p}{2})_x \!\vert^2\Big) dx\!\!+\!\!\mu_1\!\!\!\int_0^L\!\!\!\!u^p dx\nonumber \\
\!\!\!\!\!&\leq&\!\!\!\!\! -\frac{2d_1(p-1)}{p^2}\int_0^L \vert (u^\frac{p}{2})_x \vert^2 dx+\Big(\frac{(p-1)\chi^2}{2d_1}+\mu_1\Big)\int_0^L u^p dx,
\end{eqnarray}
where we have used Young's inequality in the third line of (\ref{213}).  We shall need the following estimate from P. 63 in \cite{LSU} and Corollary 1 in \cite{CKWW} with space dimension $d=1$ due to Gagliardo--Ladyzhenskaya--Nirenberg inequality that for any $u\in H^1(0,L)$ and any $\epsilon>0$, there exists some $C_0>0$ which only depends on $L$ such that
\begin{equation}\label{214}
\Vert u^\frac{p}{2}\Vert^2_{L^2(0,L)}\leq \epsilon \Vert (u^\frac{p}{2})_x \Vert^2_{L^2(0,L)}+C_0 \big(1+ \epsilon^{-\frac{1}{2}}\big)\Vert u^\frac{p}{2}\Vert^2_{L^1(0,L)}.
\end{equation}
Choosing $\epsilon=\frac{2d_1^2(p-1)}{p^2((p-1)\chi^2+2\mu_1d_1)}$ in (\ref{214}) such that $\frac{d_1(p-1)}{p^2\epsilon}=\frac{(p-1)\chi^2}{2d_1}+\mu_1$, we have that
\begin{eqnarray}\label{215}
\Big(\frac{(p-1)\chi^2}{2d_1}+\mu_1\Big)\int_0^L u^p dx
\!\!\!\!\! &=&\!\!\!\!\! \frac{2d_1 (p-1)}{p^2\epsilon}  \int_0^L u^p dx-\Big(\frac{(p-1)\chi^2}{2d_1}+\mu_1\Big)\int_0^L u^p dx\nonumber \\
\!\!\!\!\! &\leq&\!\!\!\!\! \frac{2d_1 (p-1)}{p^2} \int_0^L \vert (u^\frac{p}{2})_x \vert^2 dx+\frac{2d_1(p-1)C_0 \big(1+ \epsilon^{-\frac{1}{2}}\big)}{p^2\epsilon}\nonumber\\
&\cdot&\!\!\!\Big(\int_0^L u^\frac{p}{2} dx\Big)^2\!\!-\!\!\Big(\frac{(p-1)\chi^2}{2d_1}+\mu_1\Big)\int_0^L u^p dx.
\end{eqnarray}
In light of (\ref{215}), (\ref{213}) gives rise to
\[\frac{d}{dt}\int_0^L u^p dx\leq - \Big(\frac{p(p-1)\chi^2}{2d_1}+ p \mu_1\Big) \int_0^L u^p dx +\frac{2d_1(p-1)C_0 \big(1+ \epsilon^{-\frac{1}{2}}\big)}{p \epsilon}\Big(\int_0^L u^\frac{p}{2} dx\Big)^2.\]
On the other hand, we can choose $p_0$ large such that for all $p\geq p_0$, $\epsilon=\frac{2d_1^2(p-1)}{p^2((p-1)\chi^2+2\mu_1d_1)}$\\
$>\big(\frac{d_1}{p\chi}\big)^2$ and $\frac{1+\epsilon^{-\frac{1}{2}}}{\epsilon}<\frac{1+\frac{p\chi}{d_1}}{d_1^2/p^2\chi^2}$, then we have that
\begin{equation}\label{216}
\frac{d}{dt}\int_0^L\!\!\! u^p dx\leq -\Big(\frac{p(p-1)\chi^2}{2d_1}+ p \mu_1\Big) \int_0^L\!\!\! u^p dx +\frac{2p(p-1)C_0\chi^2(1+\frac{p\chi}{d_1})}{d_1}\Big(\int_0^L u^\frac{p}{2}dx\Big)^2.
\end{equation}
Denote $\kappa=\frac{p(p-1)\chi^2}{2d_1}+ p \mu_1$.  For each $T\in(0,\infty)$, we solve the differential inequality (\ref{216}) to have that for all $p\geq p_0$
\begin{eqnarray}\label{217}
\int_0^L u^p dx   &\leq&  e^{-\kappa t} \int_0^L u_0^p dx+4C_0\Big(1+\frac{p\chi}{d_1}\Big)\int_0^t e^{-\kappa(t-s)}\Big(\int_0^L u^\frac{p}{2} dx\Big)^2ds \nonumber\\
&\leq& e^{-\kappa t}\int_0^L u_0^p dx+4C_0\left(1+\frac{p\chi}{d_1}\right) \sup_{t\in(0,T)}\Big(\int_0^L u^\frac{p}{2}dx\Big)^2. \\ \nonumber
\end{eqnarray}
We now use the Moser--Alikakos iteration \cite{A0} to establish $L^\infty$-estimate of $u$ through (\ref{217}).  To this end, we let
\[M(p)=\max\big\{\Vert u_0 \Vert_{L^\infty},~\sup_{t\in(0,T)}\Vert u(\cdot,t) \Vert_{L^p}  \big\},\]
then (\ref{217}) implies that
\[M(p)\leq  \left(C_1+\frac{C_1 p\chi}{d_1}\right)^\frac{1}{p}M(p/2),\forall p>p_0,\]
where $C_1$ is a positive constant dependent on $C_0$ and $L$.  Taking $p=2^i$, $i=1,2,...$ and choosing $p_0=2^{i_0}$ for $i_0$ large but fixed, we have
\begin{eqnarray}\label{218}
M(2^i)\!\!\!\!\! &\leq&\!\!\!\!\!\!  \Big(\!C_1\! +\!\frac{C_1 \chi 2^i}{d_1}\!\Big)^{2^{-i}}\!\!\!\!\!\!\!M(2^{i-1})\!\!\leq\!\! \Big(\!C_1 \!+\!\frac{C_1 \chi 2^i}{d_1}\!\Big)^{2^{-i}}\!\!\!\!\!\! \Big(\!C_1\! +\!\frac{C_1 \chi 2^{i-1}}{d_1}\!\Big)^{2^{-(i-1)}}\!\!\!\!\!\!M(2^{i-2}) \nonumber\\
\!\!\!\!\! &\leq&\!\!\!\!\! M(2^{i_0}) \prod_{j=i_0+1}^i \Big(C_1+\frac{C_1\chi 2^j}{d_1}\Big)^{2^{-j}}\!\!\!\!\!\!\!\leq M(2^{i_0}) \prod_{j=i_0+1}^i \Big(C_12^j+\frac{C_1\chi 2^j}{d_1}\Big)^{2^{-j}} \nonumber\\
\!\!\!\!\! &\leq&\!\!\!\!\!  M(2^{i_0}) \Big(C_1+\frac{C_1\chi }{d_1}\Big) \prod_{j=i_0+1}^i (2^j)^{2^{-j}} \leq M(2^{i_0}) \Big(C_1+\frac{C_1\chi }{d_1}\Big) 2^{\sum\limits_{j=i_0+1}^i j2^{-j}}  \nonumber\\
\!\!\!\!\! &\leq&\!\!\!\!\!  C\Big(1+\frac{\chi }{d_1}\Big) M(2^{i_0}),
\end{eqnarray}
where $C$ is a constant that only depends on $L$ and $M(2^{i_0})$ is bounded for all $t\in(0,\infty)$ thanks to (\ref{28}).  Sending $i\rightarrow \infty$ in (\ref{218}), we finally conclude from Lemma \ref{lemma24} and (\ref{218}) that
\begin{equation}
\Vert u(\cdot,t) \Vert_{L^\infty}\leq \lim_{i\rightarrow \infty} M(2^i)\leq CM(2^{i_0}),\forall t\in[0,\infty).
\end{equation}
By the same calculations we can show that $\Vert v(\cdot,t) \Vert_{L^\infty}$ is uniformly bounded for all $t\in(0,\infty)$.  This completes the proof of Theorem \ref{theorem25}.
\end{proof}

\section{Existence and stability of nonconstant positive steady states}\label{section3}

In this section, we consider the stationary system of (\ref{11}) in the following form
\begin{equation}\label{31}
\left\{
\begin{array}{ll}
(d_1u'-\chi u w')'+\mu_1(1-u-a_1v)u=0,&x \in (0,L), \\
(d_2v'-\xi v w')'+\mu_2(1-a_2u-v)v=0,&x \in (0,L),\\
w''-\lambda w +u+v=0,&x \in (0,L),\\
u'(x)=v'(x)=w'(x)=0,&x=0,L,
\end{array}
\right.
\end{equation}
where $'$ denotes the derivative taken with respect to $x$.  From the viewpoint of mathematical modeling, it is interesting to study the aggregation of cellular organisms through chemotaxis and interspecific competitions.  For this purpose, we study the existence and stability of nonconstant positive steady states of (\ref{31}) which can be used to model the cellular aggregations.  In the absence of chemotaxis, i.e., when $\chi=\xi=0$, (\ref{31}) becomes the classical diffusive Lotka--Volterra competition system and it has only constant nonnegative solutions according to the well--known results of K. Kishimoto and H. Weinberger \cite{KiW}.  See \cite{WGY} for detailed discussions on the classical Lotka--Volterra systems.  We are concerned with the effect of chemotaxis on the existence of nonconstant positive solutions to (\ref{31}).  Unlike random movements, a chemotactic process is anti--diffusion and it has the effect of destabilizing the spatially homogeneous solutions.  Then spatially inhomogeneous solutions may arise through bifurcation as the homogeneous solution becomes unstable.  We want to point out that advection-- or chemotaxis--driven instability has been investigated by various authors in \cite{MOW}, \cite{WGY}, etc.

\subsection{Linearized stability analysis of the homogeneous solution $(\bar u,\bar v,\bar w)$}

To study the mechanism through which spatially inhomogeneous solutions of (\ref{31}) emerge, we carry out the standard linearized stability of $(\bar u,\bar v,\bar w)$ given by (\ref{13}), viewed as an equilibrium of (\ref{11}).  To this end, we take $(u,v,w)=(\bar u,\bar v,\bar w)+(U,V,W)$, where $U$, $V$ and $W$ are small perturbations from $(\bar u,\bar v,\bar w)$, then $(U,V,W)$ satisfies the following system
\begin{equation*}
\left\{
\begin{array}{ll}
U_t\approx (d_1U'-\chi\bar u W')'-\mu_1\bar uU-\mu_1a_1\bar uV,&x \in (0,L),t>0      \\
V_t\approx (d_2V'-\xi\bar v W')'-\mu_2a_2\bar vU-\mu_2\bar vV,&x \in (0,L),t>0      \\
W_t\approx W''-\lambda W+U+V,&x \in (0,L),t>0      \\
U'(x)=V'(x)=W'(x)=0,&x=0,L,t>0.
\end{array}
\right.
\end{equation*}
According to the standard linearized stability analysis--see \cite{Si} for example, we know that the stability of $(\bar u,\bar v,\bar w)$ can be determined by the following matrix,
\begin{equation}\label{32}
\left(
  \begin{array}{ccc}
    -d_1\Lambda-\mu_1\bar u & -\mu_1a_1\bar u & \chi\bar u\Lambda\\
    -\mu_2a_2\bar v & -d_2\Lambda-\mu_2\bar v & \xi\bar v\Lambda\\
     1& 1 & -\Lambda-\lambda
\end{array}
\right),
\end{equation}
where $\Lambda=(\frac{k\pi}{L})^2>0,k=1,2,...$, are the eigenvalues of $-\frac{d^2}{dx^2}$ on $(0,L)$ under the Neumann boundary condition.  We have the following result on the linearized instability of $(\bar u,\bar v,\bar w)$.
\begin{proposition}\label{proposition1}
Assume that condition (\ref{14}) is satisfied, the positive constant solution $(\bar u,\bar v,\bar w)$ is unstable with respect to (\ref{11}) if
\begin{equation}\label{33}
\chi\geq\chi_0=\min_{k\in \mathbb{N}^+}\{\tilde \chi_k,\hat \chi_k\},
\end{equation}
and it is locally asymptotically stable if $\chi<\chi_0=\min_{k\in \mathbb{N}^+}\{\tilde \chi_k,\hat \chi_k\}$,
where
\begin{align}\label{34}
\tilde \chi_k=&\frac{\big( \big(d_1(\frac{k\pi}{L})^2+\mu_1\bar u\big) \big(d_2(\frac{k\pi}{L})^2+\mu_2\bar v\big)-a_1a_2\mu_1\mu_2\bar u\bar v  \big) \big((\frac{k\pi}{L})^2+\lambda \big)}{d_2\bar u(\frac{k\pi}{L})^4+(1-a_2)\mu_2\bar u\bar v(\frac{k\pi}{L})^2}\nonumber \\
&-\frac{\xi \big(d_1\bar v(\frac{k\pi}{L})^4+(1-a_1)\mu_1\bar u\bar v(\frac{k\pi}{L})^2\big)}{d_2\bar u(\frac{k\pi}{L})^4+(1-a_2)\mu_2\bar u\bar v(\frac{k\pi}{L})^2},
\end{align}
and
\begin{equation}\label{35}
\hat \chi_k=\frac{-\xi G_2+H_1{H_2}^2+{H_1}^2H_2+H_2H_3}{G_1},
\end{equation}
where
\[G_1=(d_1+1)\bar u(\frac{k\pi}{L})^4+(\lambda+\mu_1\bar u+a_2\mu_2\bar v)\bar u(\frac{k\pi}{L})^2 ,\]
\[G_2=(d_2+1)\bar v(\frac{k\pi}{L})^4+(\lambda+a_1\mu_1\bar u+\mu_2\bar v)\bar v(\frac{k\pi}{L})^2,\]
\[H_1=(\frac{k\pi}{L})^2+\lambda,~H_2=(d_1+d_2)(\frac{k\pi}{L})^2+\mu_1\bar u+\mu_2\bar v,\]
and
\[H_3=\Big(d_1(\frac{k\pi}{L})^2+\mu_1\bar u\Big) \Big(d_2(\frac{k\pi}{L})^2+\mu_2\bar v \Big)-a_1a_2\mu_1\mu_2\bar u\bar v.\]
\end{proposition}
\begin{proof}
By the principle of the linearized stability--see Theorem 5.2 in \cite{Si}, $(\bar u,\bar v,\bar w)$ is asymptotically stable with respect to (\ref{11}) if and only if the real parts of all eigenvalues of the matrix (\ref{32}) are negative.  To this end, we first see that its characteristic polynomial reads
\begin{equation*}
\sigma^3+\alpha_2(k)\sigma^2+\alpha_1(k)\sigma+\alpha_0(k)=0,
\end{equation*}
where
\begin{equation*}
\alpha_2(k)= (d_1+d_2+1) \Big(\frac{k\pi}{L}\Big)^2+\mu_1\bar u+\mu_2\bar v+\lambda,
\end{equation*}
\begin{align*}
\alpha_1(k)=&\Big(\Big(\frac{k\pi}{L}\Big)^2+\lambda \Big) \Big(\Big(d_1+d_2\Big)\Big(\frac{k\pi}{L}\Big)^2
+\mu_1\bar u+\mu_2\bar v \Big)-a_1a_2\mu_1\mu_2\bar u\bar v \nonumber \\
&-(\chi\bar u+\xi\bar v)\Big(\frac{k\pi}{L}\Big)^2+\Big(d_1\Big(\frac{k\pi}{L}\Big)^2+\mu_1\bar u\Big) \Big(d_2\Big(\frac{k\pi}{L}\Big)^2+\mu_2 \bar v\Big),
\end{align*}
and
\begin{align*}
\alpha_0(k)=\!&-\chi\bar u\Big(\frac{k\pi}{L}\Big)^2\!\Big(d_2\Big(\frac{k\pi}{L}\Big)^2
\!\!+\!(1\!-\!a_2)\mu_2\bar v \Big)\!-\!\xi\bar v\Big(\frac{k\pi}{L}\Big)^2 \!\Big(d_1\Big(\frac{k\pi}{L}\Big)^2\!+\!(1\!-\!a_1)\mu_1\bar u \Big) \nonumber\\
&-\!a_1a_2\mu_1\mu_2\bar u\bar v\Big(\!\Big(\frac{k\pi}{L}\Big)^2\!\!\! +\!\lambda \Big)\!+\!\Big(d_1\Big(\frac{k\pi}{L}\Big)^2\!\!+\!\mu_1\bar u\Big)\Big(\!d_2\Big(\frac{k\pi}{L}\Big)^2\!\!+\!\mu_2\bar v\Big)\!\Big(\Big(\!\frac{k\pi}{L}\Big)^2\!\! +\!\lambda \Big).
\end{align*}
We see that $\alpha_2(k)>0$ since $d_i$, $\mu_i$, $i=1,2$ and $\lambda$ are all positive.  According to the Routh--Hurwitz conditions, or Corollary 2.2 in \cite{LSW}, the real parts of all eigenvalues $\sigma_k$ to (\ref{32}) are negative hence $(\bar u,\bar v,\bar w)$ is locally stable if and only if for all $k\in \mathbb N^+$
\[\alpha_0(k)>0, \alpha_1(k)>0, \text{~and~}\alpha_1(k)\alpha_2(k)-\alpha_0(k)>0,\]
while the real parts of some eigenvalues are nonnegative if one of the conditions above fails for some $k\in \mathbb N^+$.  Moreover, we can show by simple calculations that $\alpha_1(k)>0$ is necessary whenever $\alpha_0(k)>0$ and $\alpha_1(k)\alpha_2(k)-\alpha_0(k)>0$.  Therefore, $(\bar u,\bar v,\bar w)$ is unstable if there exists some $k\in \mathbb N^+$ such that one of the following conditions is satisfied,
\[\alpha_0(k)<0,\text{ or }\alpha_1(k)\alpha_2(k)-\alpha_0(k)<0.\]
We have from simple calculations that $\chi>\tilde \chi_k$ if and only if $\alpha_0(k)<0$ and $\alpha_1(k)\alpha_2(k)-\alpha_0(k)<0$ if and only if $\chi>\hat \chi_k$.  Therefore $(\bar u,\bar v,\bar w)$ is unstable if $\chi>\chi_0$ according to (S1) of Corollary 2.2 in \cite{LSW}.  Similarly, we can show that $(\bar u,\bar v,\bar w)$ is locally asymptotically stable if $\chi<\chi_0$.
\end{proof}

Proposition \ref{proposition1} states that the homogeneous steady state $(\bar u,\bar v,\bar w)$ loses its stability at $\chi=\chi_0$.  We shall see in our coming analysis that the stability is lost to stable bifurcating solutions in most cases.  It is easy to see that large $\xi$ can also destabilize the homogeneous solution and we investigate the effect of $\chi$ in this paper with $\xi$ being fixed.  It is worthwhile to mention that Proposition 1 also holds for multi-dimensional domain $\Omega\subseteq\mathbb{R}^N$, $N\geq 2$, with $(\frac{k\pi}{L})^2$ being replaced by the $k$--eigenvalue of $-\Delta $ under the Neumann boundary condition.

\begin{remark}\label{remark1}
If $\chi_0=\tilde \chi_k$, then $\alpha_0(k)=0$ and the eigenvalues of (\ref{32}) are $\tilde \sigma_1(k)=0$, $\tilde \sigma_2(k)=\frac{-\alpha_2(k)+\sqrt{\alpha^2_2(k)-4\alpha_1(k)}}{2}$ and $\tilde \sigma_3(k)=\frac{-\alpha_2(k)-\sqrt{\alpha^2_2(k)-4\alpha_1(k)}}{2}$; if $\chi_0=\hat \chi_k$, then $\alpha_0(k)=\alpha_1(k)\alpha_2(k)$ and the eigenvalues of (\ref{32}) are $\hat \sigma_1(k)=-\alpha_2(k)<0$, $\hat \sigma_2(k)=\sqrt{-\alpha_1(k)}$ and $\hat \sigma_3(k)=-\sqrt{-\alpha_1(k)}$.  We shall see in our coming analysis and numerical simulations that $(\bar u,\bar v,\bar w)$ loses its stability through either Hopf bifurcation or steady state bifurcation.  This depends on whether $\chi_0=\min_{k\in \mathbb{N}^+}\{\tilde \chi_k,\hat \chi_k\}$ in (\ref{33}) is achieved at $\min_{k\in \mathbb{N}^+} \tilde \chi_k$ or $\min_{k\in \mathbb{N}^+} \hat \chi_k$.  We divide our discussions into the following cases:

\begin{itemize}
\item Case 1. $\chi_0=\tilde \chi_{k_0}<\min_{k\in\mathbb N^+} \hat\chi_k$.  Since $(\bar u,\bar v,\bar w)$ is unstable for all $\chi>\chi_0$, we know that (\ref{32}) has at least one eigenvalue $\tilde \sigma_i$ with positive real part (in particular) when $\chi=\tilde \chi_k$, $k\neq k_0$.  This fact is very important in our proof of the stability of nonconstant positive solutions to (\ref{31}) that bifurcate from $(\bar u,\bar v,\bar w,\tilde \chi_k)$.  In particular, it implies that the only stable bifurcating solutions must be on the branch around $(\bar u,\bar v,\bar w,\tilde \chi_{k_0})$ that turns to the right and all the rest branches are always unstable.  See Proposition \ref{proposition2} or Theorem \ref{theorem32} in Section 3.3.
\item Case 2. $\chi_0=\hat \chi_{k_1}<\min_{k\in\mathbb N^+}\tilde \chi_k$.  Similar as in Case 1, we can show that the bifurcating solutions around $(\bar u,\bar v,\bar w,\tilde \chi_k)$ are unstable for all $k\in\mathbb N^+$.  We claim that $\alpha_1(k_1)\geq0$ with $\chi=\hat \chi_{k_1}$, therefore (\ref{32}) has three eigenvalues $\hat\sigma_1=-\alpha_2(k_1)$, $\hat\sigma_2=\sqrt{\alpha_1(k_1)}i$ and $\hat\sigma_3=-\sqrt{\alpha_1(k_1)}i$.  This indicates the possibility of a Hopf bifurcation hence the emergence of spatial time--periodic patterns in (\ref{11}) when $\chi=\hat \chi_{k_1}$.  To prove this claim, we argue by contradiction and assume that $\alpha_1(k_1)<0$, therefore $\hat \sigma_2(k_1)=\sqrt{-\alpha_1(k_1)}>0$ for $\chi=\hat \chi_{k_1}$ and $\hat \sigma_2>0$ if $\chi$ is slightly smaller than $\chi_{k_1}$, however this indicates that $(\bar u,\bar v,\bar w)$ is unstable for $\chi<\hat \chi_{k_1}=\chi_0$ which is a contradiction.  Our numerical simulations illustrate the existence of spatial--temporal periodic patterns in (\ref{11})--see Figure \ref{fig4} and Figure \ref{fig5}.  However, rigorous Hopf bifurcation analysis of (\ref{11}) is out of the scope of this paper and we postpone it for further studies.
\item Case 3. $\chi_0=\min_{k\in\mathbb N^+}\tilde \chi_k=\min_{k\in\mathbb N^+}\hat \chi_k$.  In this case, matrix (\ref{32}) has three eigenvalues $\sigma_1=-\alpha_2(k)<0$, $\sigma_2=\sigma_3=0$, and the linear stability of $(\bar u,\bar v,\bar w)$ is lost for an zero eigenvalue with multiplicity 2.  This is also true if $\tilde \chi_k=\hat \chi_k$ $k\in\mathbb N^+$.  This inhibits the application of our steady state bifurcation analysis which requires the null--space of (\ref{32}) to be one--dimensional.  Therefore we assume that $\tilde \chi_k\neq\hat \chi_k$, $\forall k\in\mathbb N^+$ in our steady state bifurcation analysis.  See Theorem \ref{theorem31}.
\end{itemize}

\end{remark}

\begin{remark}\label{remark2}
The critical value $\chi_0=\chi_0(\xi)$ decreases as $\xi$ increases;  moreover $\chi_0<0$ if $\xi$ is sufficiently large.  Therefore only one of $\chi$ and $\xi$ is needed to be large to destabilize $(\bar u,\bar v,\bar w)$.  If $\xi<0$, i.e., species $v$ is repulsive to the chemical gradient, then $\chi$ needs to be large to destabilize $(\bar u,\bar v,\bar w)$.  Moreover, the local stability analysis suggests that chemo--attraction destabilizes constant steady states and the chemo--repulsion stabilizes constant steady states.  The constant solution is always stable when $\chi<0$ and $\xi<0$.  Therefore, we surmise that $(\bar u,\bar v,\bar w)$ is also a global attractor of (\ref{11}) if $\chi<0$ and $\xi<0$, though the positiveness of $\chi$ and $\xi$ is required in \cite{TW}.
\end{remark}

\subsection{Bifurcation of steady states}
We now proceed to perform rigorous steady state bifurcation analysis to seek nonconstant positive solutions to
the stationary system (\ref{31}).  As we have shown above, the chemo--attraction rates $\chi$ and $\xi$ have the effect of destabilizing $(\bar u,\bar v,\bar w)$ which becomes unstable if $\chi$ surpasses $\chi_0$.  The linearized instability of $(\bar u,\bar v,\bar w)$ in (\ref{11}) is insufficient to guarantee the existence of spatially inhomogeneous steady states.  Therefore we are concerned with the conditions under which the constant solution $(\bar u,\bar v,\bar w)$ loses its stability to spatially inhomogeneous solutions through steady state bifurcation.  Clearly, the emergence of spatially inhomogeneous solutions is due to the effect of large chemo--attraction rate $\chi$ and we refer these as advection--induced patterns in the sense of Turing's instability.

We shall apply the bifurcation theory of Crandall--Rabinowitz \cite{CR} to seek nonconstant positive solutions to (\ref{31}).  To this end, we introduce the Hilbert space
\[\mathcal{X}=\{w\in H^2(0,L)\mid w'(0)=w'(L)=0\},\]
and convert (\ref{31}) into the following abstract form
$$\mathcal{F}(u,v,w,\chi)=0, \ (u,v,w,\chi)\in \mathcal{X}\times \mathcal{X}\times \mathcal{X}\times \mathbb{R},$$
where
\begin{equation*}
\mathcal{F}(u,v,w,\chi)=\left(
  \begin{array}{c}
(d_1u'-\chi u w')'+\mu_1(1-u-a_1v)u \\
(d_2v'-\xi v w')'+\mu_2(1-a_2u-v)v \\
 w''-\lambda w +u+v
\end{array}
\right).
\end{equation*}
It is easy to see that $\mathcal{F}(\bar u,\bar v,\bar w,\chi)=0$ for any $\chi \in \mathbb{R}$ and $\mathcal{F}:\mathcal{X} \times \mathcal{X}\times \mathcal{X}\times \mathbb{R} \rightarrow \mathcal{Y}\times \mathcal{Y}\times \mathcal{Y}$ is analytic for $\mathcal{Y}=L^2(0,L)$; moreover, it follows from straight calculations that for any fixed $(\hat{u},\hat{v},\hat{w})\in \mathcal{X}\times \mathcal{X}\times \mathcal{X}$, the Fr$\acute{\text{e}}$chet derivative of $\mathcal{F}$ is given by
\begin{align}\label{36}
&D_{(u,v,w)}\mathcal{F}(\hat{u},\hat{v},\hat{w},\chi)(u,v,w) \nonumber \\
= & \left(
  \begin{array}{c}
    \big(d_1u'-\chi( u\hat{w}'+\hat{u}w')\big)'+\mu_1(1-2 \hat{u}-a_1 \hat{v})u-a_1\mu_1\hat{u}v \\
    \big(d_2v'-\xi( v\hat{w}'+\hat{v}w')\big)'+\mu_2(1-a_2 \hat{u}-2 \hat{v})v-a_2\mu_2 u\hat{v}  \\
    w''-\lambda w+u+v
\end{array}
\right).
\end{align}
We claim that the derivative $D_{(u,v,w)}\mathcal{F}(\hat{u},\hat{v},\hat{w},\chi):\mathcal{X}\times \mathcal{X}\times \mathcal{X}\times \mathbb{R}\rightarrow \mathcal{Y} \times \mathcal{Y}\times \mathcal{Y}$ is a Fredholm operator with zero index.  To see this, we denote $\textbf{u}=(u,v,w)^\text{T}$ and write (\ref{36}) into
\[D_{(u,v,w)}\mathcal{F}(\hat{u},\hat{v},\hat{w},\chi)(u,v,w)=\mathcal{A}_0(\textbf{u}) \textbf{u}''+\textbf{F}_0(x,\textbf{u},  \textbf{u}'),\]
where
\[\mathcal{\textbf{A}}_0(\textbf{u})=\begin{pmatrix}
d_1 & 0 & -\chi \hat u  \\
0 &  d_2 &-\xi \hat v \\
0 & 0 & 1
\end{pmatrix}\]
and
\[\mathcal{\textbf{F}}_0=\begin{pmatrix}
-\chi(u\hat{w}')'-\chi \hat{u}'w'+\mu_1(1-2 \hat{u}-a_1 \hat{v})u-a_1\mu_1\hat{u}v \\
-\xi(v\hat{w}')'-\xi \hat{v}'w'+\mu_2(1-a_2 \hat{u}-2 \hat{v})v-a_2\mu_2\hat{v}u \\
-\lambda w +u+v
\end{pmatrix}.\]  It is obvious that operator (\ref{36}) is strictly elliptic since all eigenvalues of $\mathcal{A}_0$ are positive.  According to Remark 2.5 (case 2) in Shi and Wang \cite{SW} with $N=1$, $\mathcal{A}_0$ satisfies the Agmon's condition, which in loose terms is equivalent as the normal ellipticity (see p.21 of \cite{A0} and Definition 2.4 in \cite{SW} for the Agmon's condition).  Therefore, $D_{(u,v,w)}\mathcal{F}(\hat u,\hat v,\hat w,\chi)$ is a Fredholm operator with zero index thanks to Theorem 3.3 and Remark 3.4 of \cite{SW}.

For local bifurcations to occur at $(\bar u,\bar v,\bar w,\chi)$, we need that the Implicit Function Theorem to fails on $\mathcal{F}$ or equivalently the following non--triviality condition is satisfied
\[\mathcal{N}\big(D_{(u,v,w)}\mathcal{F}(\bar u,\bar v,\bar w,\chi)\big)\neq \{0\},\]
where $\mathcal{N}$ denotes the null space.  Taking $(\hat{u},\hat{v},\hat{w})=(\bar u,\bar v,\bar w)$ in (\ref{36}), we have that
\begin{equation}\label{37}
D_{(u,v,w)}\mathcal{F}(\bar u,\bar v,\bar w,\chi)(u,v,w)=
  \left(
  \begin{array}{c}
    d_1u''-\chi \bar uw''-\mu_1\bar uu-a_1\mu_1\bar uv \\
    d_2v''-\xi \bar vw''-\mu_2\bar vv-a_2\mu_2\bar vu \\
    w''-\lambda w+u+v
\end{array}
\right),
\end{equation}
and its null space consists of solutions to the following system
\begin{eqnarray*}
\left\{
\begin{array}{ll}
d_1u''-\chi\bar uw''-\mu_1\bar uu-a_1\mu_1\bar uv=0,&x\in(0,L),\\
d_2v''-\xi\bar vw''-\mu_2\bar vv-a_2\mu_2\bar vu=0,&x\in(0,L), \\
w''-\lambda w+u+v=0,&x\in(0,L),\\
u'(x)=v'(x)=w'(x)=0,&x=0,L.
\end{array}
\right.
\end{eqnarray*}
Substituting the eigen--expansions
\[u(x)=\sum_{k=0}^{\infty}t_k\cos\frac{k\pi x}{L}, \ v(x)=\sum_{k=0}^{\infty}s_k\cos\frac{k\pi x}{L},\ w(x)=\sum_{k=0}^{\infty}r_k\cos\frac{k\pi x}{L},\]
into the system above, we have that $(t_k,s_k,r_k)$ satisfies
\begin{equation}\label{38}
\left(
\begin{array}{ccc}
    -(\!d_1(\!\frac{k\pi}{L}\!)^2\!+\!\mu_1\bar u\!) &\!\! -a_1\mu_1\bar u &\!\! \chi\bar u(\!\frac{k\pi}{L}\!)^2\\
    -a_2\mu_2\bar v &\!\! -(\!d_2(\!\frac{k\pi}{L}\!)^2\!+\!\mu_2\bar v\!) &\!\! \xi\bar v(\!\frac{k\pi}{L}\!)^2\\
     1&\!\! 1 &\!\! -(\!(\!\frac{k\pi}{L}\!)^2+\lambda\!)
\end{array}
\right)
\!\!\left(
\begin{array}{c}t_k\\
s_k\\
r_k
\end{array}
\!\!\right)\!\!=\!\!
\!\!\left(
\begin{array}{c}0\\
0\\
0
\end{array}
\!\!\right).
\end{equation}
First of all, $k = 0$ can be easily ruled out in (\ref{38}) if $(t_k,s_k,r_k)$ is nonzero.  For $k\in \mathbb{N}^+$, the coefficient matrix of (\ref{38}) is singular and it follows from simple calculations that
\begin{align}\label{39}
\chi=\chi_k=&\frac{ \big(\big(d_1(\frac{k\pi}{L})^2+\mu_1\bar u\big)\big(d_2(\frac{k\pi}{L})^2+\mu_2\bar v\big)-a_1a_2\mu_1\mu_2\bar u\bar v \big)\big((\frac{k\pi}{L})^2+\lambda\big)}{d_2\bar u(\frac{k\pi}{L})^4+(1-a_2)\mu_2\bar u\bar v(\frac{k\pi}{L})^2} \nonumber \\
&-\frac{\xi \big(d_1\bar v(\frac{k\pi}{L})^4+(1-a_1)\mu_1\bar u\bar v(\frac{k\pi}{L})^2\big)}{d_2\bar u(\frac{k\pi}{L})^4+(1-a_2)\mu_2\bar u\bar v(\frac{k\pi}{L})^2}, k\in \mathbb{N}^+.
\end{align}
It is easy to see that $\chi_k>0$ if and only if
\begin{equation*}
\xi<\frac{\big(\big(d_1(\frac{k\pi}{L})^2+\mu_1\bar u\big)\big(d_2(\frac{k\pi}{L})^2+\mu_2\bar v\big)-a_1a_2\mu_1\mu_2\bar u\bar v\big)\big((\frac{k\pi}{L})^2+\lambda\big)}{d_1\bar v(\frac{k\pi}{L})^4+(1-a_1)\mu_1\bar u\bar v(\frac{k\pi}{L})^2}.
\end{equation*}
We notice that $\chi_k$ is the same as $\tilde \chi_k$ given in (\ref{34}).  Moreover, according to our discussions in Case 3 of Remark \ref{remark1}, if $\chi_k\neq \hat\chi_k$, the null space $D_{(u,v,w)}\mathcal{F}(\bar u,\bar v,\bar w,\chi)$ has dimension 1 and it is spanned by $\{(\bar u_k,\bar v_k,\bar w_k)\}_{k=1}^\infty$, where for each $k\in \mathbb{N}^+$
\begin{equation}\label{310}
(\bar u_k,\bar v_k,\bar w_k)=(P_k,Q_k,1)\cos\frac{k\pi x}{L},
\end{equation}
with
\begin{equation}\label{311}
P_k=\frac{\big(d_2(\frac{k\pi}{L})^2+\mu_2\bar v\big)\big((\frac{k\pi}{L})^2+\lambda\big)-\xi\bar v(\frac{k\pi}{L})^2}{d_2(\frac{k\pi}{L})^2+(1-a_2)\mu_2\bar v},
\end{equation}
and
\begin{equation}\label{312}
Q_k=\frac{\xi\bar v(\frac{k\pi}{L})^2-a_2\mu_2\bar v\big((\frac{k\pi}{L})^2+\lambda\big)}{d_2(\frac{k\pi}{L})^2+(1-a_2)\mu_2\bar v}.
\end{equation}

We want to point out that the bifurcation value $\chi_k$ in (\ref{39}) is greater than or equal to $\chi_0$ given by (\ref{33}) in Proposition \ref{proposition1} at which $(\bar u,\bar v,\bar w)$ loses its stability.  If $\chi_0=\min_{k\in \mathbb{N}^+}\chi_k<\min_{k\in \mathbb{N}^+}\hat \chi_k$, then we shall see that $(\bar u,\bar v,\bar w)$ loses its stability to spatially inhomogeneous solutions that emerge through steady state bifurcations.  However, if $\chi_0=\min_{k\in \mathbb{N}^+} \hat \chi_k<\min_{k\in \mathbb{N}^+} \chi_k$, it is unclear to us what happens in this case.  The discussions in Case 3 of Remark \ref{remark1} and our numerical simulations in Section \ref{section4} suggest that a Hopf bifurcation may occur at $(\bar u,\bar v,\bar w)$ to create spatial time--periodic patterns.  However, the rigorous analysis is out of the scope of this paper and we refer the reader to the discussions in Section 2 of \cite{LSW}.

Having the potential bifurcation value $\chi_k$ in (\ref{39}), we now verify that local bifurcation does occur at $(\bar{u},\bar{v},\bar{w},\chi_k)$ for each $k\in \mathbb{N}^+$ in the following theorem, which establishes the existence of nonconstant positive solutions to (\ref{31}).
\begin{theorem}\label{theorem31}
Let $d_i$, $\mu_i$, $i=1,2$ and $\lambda$ be positive constants and assume that $0<a_1, a_2<1$.  Suppose that for all positive different integers $k,j\in \mathbb N^+$,
\begin{equation}\label{313}
 \chi_k\neq \chi_j, \forall k\neq j\text{ and }\chi_k\neq \hat \chi_k, \forall k\in \mathbb N^+,
\end{equation}
where $\hat \chi_k$ and $\chi_k$ are defined in (\ref{35}) and (\ref{39}) respectively, and $(\bar u,\bar v,\bar w)$ is the positive equilibrium given in (\ref{13}).  Then for each $k\in \mathbb N^+$, there exists a branch of spatially inhomogeneous solutions $(u,v,w)\in \mathcal{X}\times \mathcal{X}\times \mathcal{X}$ of (\ref{31}) that bifurcate from $(\bar{u},\bar{v},\bar{w})$ at $\chi=\chi_k$.  Moreover, for $\delta>0$ being small, the bifurcating branch $\Gamma_k(s)$ around $(\bar{u},\bar{v},\bar{w},\chi_k)$ can be parameterized as
\begin{equation}\label{314}
\chi_k(s)=\chi_k+O(s),s\in(-\delta,\delta)
\end{equation}
and
\begin{equation}\label{315}
(u_k(s,x),v_k(s,x),w_k(s,x))=(\bar{u},\bar{v},\bar{w})+s(\bar u_k,\bar v_k,\bar w_k) +o(s),s\in(-\delta,\delta),
\end{equation}
where $(\bar u_k,\bar v_k,\bar w_k)$ is given by (\ref{310}) and
$(u_k(s,x),v_k(s,x),w_k(s,x))-(\bar{u},\bar{v},\bar{w})-s(\bar u_k,\bar v_k,\bar w_k)$ is in the closed complement $\mathcal{Z}$ of $\mathcal{N}\big(D_{(u,v,w)}\mathcal{F}(\bar u,\bar v,\bar w,\chi)\big)$ defined as
\begin{equation}\label{316}
\mathcal{Z}=\big\{(u,v,w)\in \mathcal{X} \times \mathcal{X}\times \mathcal{X}  \big \vert \int_0^L u \bar u_k+v\bar v_k+w\bar w_k dx=0\big\};
\end{equation}
furthermore, $(u_k(s,x),v_k(s,x),w_k(s,x),\chi_k(s))$ solves system (\ref{31}) and all nontrivial solutions of (\ref{31}) near the bifurcation point ($\bar{u},\bar{v},\bar{w},\chi_k)$ must stay on the curve $\Gamma_k(s)$, $s\in(-\delta,\delta)$.
\end{theorem}
\begin{proof}
To apply the local theory of M. Crandall and P. Rabinowitz in \cite{CR}, we now only need to verify the following transversality condition
\begin{equation}\label{317}
\frac{d}{d \chi} \left(D_{(u,v,w)}\mathcal{F}(\bar{u},\bar{v},\bar{w},\chi)\right)(\bar u_k,\bar v_k,\bar w_k)\vert_{\chi=\chi_k} \notin \mathcal{R}(D_{(u,v,w)}\mathcal{F}(\bar{u},\bar{v},\bar w,\chi_k)),
\end{equation}
where $\mathcal{R}(\cdot)$ denotes the range of the operator.  If (\ref{317}) fails, there must exist a nontrivial solution $(u,v,w)$ to the following problem
\begin{eqnarray}\label{318}
\left\{
\begin{array}{ll}
d_1u''-\chi_k\bar uw''-\mu_1\bar uu-\mu_1a_1\bar uv=\big(\frac{k\pi}{L}\big)^2\bar u\cos \frac{k\pi x}{L},&x\in(0,L),\\
d_2v''-\xi\bar vw''-\mu_2\bar vv-\mu_2a_2\bar vu=0,&x\in(0,L), \\
w''-\lambda w+u+v=0,&x\in(0,L),\\
u'(x)=v'(x)=w'(x)=0,&x=0,L.
\end{array}
\right.
\end{eqnarray}
Multiplying both equations in (\ref{318}) by $\cos \frac{k\pi x}{L}$ and integrating them over $(0,L)$ by parts, we obtain
\begin{equation}\label{319}
\!\!\begin{pmatrix}
\!\!-(d_1(\!\frac{k\pi}{L}\!)^2+\mu_1\bar u\!) & \!-a_1\mu_1\bar u & \!\chi_k\bar u(\frac{k\pi}{L})^2 \\
-a_2\mu_2\bar v & \!-(d_2(\frac{k\pi}{L})^2+\mu_2\bar v)\! & \xi\bar v(\frac{k\pi}{L})^2\\
1& \!1\! & -((\frac{k\pi}{L})^2+\lambda)
\end{pmatrix}\!\!
\!\begin{pmatrix}
\!\int_0^L u \cos \frac{k \pi x}{L} dx\\
\!\int_0^L v \cos \frac{k \pi x}{L} dx\\
\!\int_0^L w \cos \frac{k \pi x}{L} dx
\end{pmatrix}\!\!=\!\!\begin{pmatrix}
\!\frac{(k\pi)^2 \bar u }{2L}\\
0\\
0
\end{pmatrix}\!\!.
\end{equation}
The coefficient matrix in (\ref{319}) is singular because of (\ref{39}), therefore we reach a contradiction and this verifies (\ref{317}).  The statements in Theorem \ref{theorem31} follow from Theorem 1.7 of \cite{CR}.
\end{proof}
We want to mention that the condition (\ref{313}) is assumed such that $(\bar u,\bar v,\bar w)$ loses its stability at $\chi_k$ for an zero eigenvalue with multiplicity one.  What happens if one of the conditions fails is out of the scope of this paper.  According to the global bifurcation theory of P. Rabinowitz \cite{Ra} and its developed version of J. Shi and X. Wang \cite{SW}, the global continuum of $\Gamma_k(s)$ will, in loose terms, either intersect with the $\chi$--axis at another bifurcating point, or extend to infinity, or meet with a \emph{singular} point.  It is important to investigate the global bifurcation of (\ref{31}) for qualitative properties of its positive solutions.  See Section \ref{section5} for more discussions.

\subsection{Stability of the bifurcating solutions}
Now we proceed to study the stability of $\big(u_k(s,x), v_k(s,x), w_k(s,x)\big)$ established in Theorem \ref{theorem31}.  Here the stability or instability is that of the bifurcating solution viewed as an equilibrium of system (\ref{11}).  To this end, we first determine the turning direction of each bifurcation branch $\Gamma_k(s)$ around $(\bar u,\bar v,\bar w,\chi_k)$, $k\in \mathbb N^{+}$.  It is easy to see that $\mathcal{F}$ is $C^4$--smooth, therefore according to Theorem 1.18 from \cite{CR}, the bifurcating solution $(u_k(s,k),v_k(s,k), w_k(s,k),\chi_k(s))$ is $C^3$--smooth of $s$ and we can write the following expansions
\begin{equation}\label{320}
\left\{
\begin{array}{ll}
u_k(s,x)=\bar u+sP_k\cos\frac{k\pi x}{L}+s^2 \varphi_1(x)+s^3\varphi_2(x)+o(s^3),\\
v_k(s,x)=\bar v+sQ_k\cos\frac{k\pi x}{L}+s^2 \psi_1(x)+s^3\psi_2(x)+o(s^3),\\
w_k(s,x)=\bar w+s \cos\frac{k\pi x}{L}+s^2 \gamma_1(x)+s^3\gamma_2(x)+o(s^3),\\
\chi_k(s)=\chi_k+s\mathcal{K}_1+s^2\mathcal{K}_2+o(s^2),
\end{array}
\right.
\end{equation}
where for $i=1,2$, $(\varphi_i, \psi_i, \gamma_i)$ belong to $\mathcal{Z}$ defined in (\ref{316}) and $\mathcal{K}_i$ are constants to be determined.  The $o(s^3)$--terms in (\ref{320}) are taken with respect to the $\mathcal{X}$--topology.

Before performing the stability analysis, we claim that the bifurcation is of pitch--fork type by showing $\mathcal{K}_1=0$.  Substituting (\ref{320}) into (\ref{31}) and equating the $s^2$--terms, we collect the following system
\begin{equation}\label{321}
\left\{
\begin{array}{ll}
d_1\varphi_{1}^{''}-\chi_k \bar u\gamma_{1}^{''}=-\chi_k P_k(\frac{k\pi}{L})^2\cos\frac{2k\pi x}{L}-\mathcal{K}_1\bar u(\frac{k\pi}{L})^2\cos\frac{k\pi x}{L}\\
\hspace{1.05in}+\mu_1 \big( \bar u(\varphi_1+a_1\psi_1) +P_k(P_k+a_1Q_k)\cos^2\frac{k\pi x}{L} \big),&x\in(0,L),\\
d_2\psi_{1}^{''}-\xi \bar v\gamma_{1}^{''}=-\xi Q_k(\frac{k\pi}{L})^2\cos\frac{2k\pi x}{L}+\mu_2 \big(\bar v(a_2\varphi_1+\psi_1)\\
\hspace{0.95in}+Q_k(Q_k+a_2P_k)\cos^2\frac{k\pi x}{L} \big),&x\in(0,L),\\
\gamma_{1}^{''}-\lambda\gamma_1+\varphi_1+\psi_1 =0,&x\in(0,L),\\
\varphi_{1}^{'}(x)=\psi_{1}^{'}(x)=\gamma_{1}^{'}(x)=0,&x=0,L.
\end{array}
\right.
\end{equation}
Multiplying the first three equations of (\ref{321}) by $\cos\frac{k\pi x}{L}$ and then integrating them over $(0,L)$ by parts, we have from simple calculations that
\begin{align}\label{322}
\frac{\bar u(k\pi)^2}{2L}\mathcal{K}_1=&\Big(d_1\Big(\frac{k\pi}{L}\Big)^2+\mu_1\bar u \Big)\int_0^L\varphi_1\cos\frac{k\pi x}{L}dx+a_1 \mu_1\bar u\int_0^L\psi_1\cos\frac{k\pi x}{L}dx \nonumber \\
&-\chi \bar u (\frac{k\pi}{L})^2\int_0^L\gamma_1\cos\frac{k\pi x}{L}dx,
\end{align}
\begin{align}\label{323}
&a_2 \mu_2\bar v\int_0^L\varphi_1\cos\frac{k\pi x}{L}dx\nonumber+\Big(d_2\Big(\frac{k\pi}{L}\Big)^2+\mu_2\bar v \Big)\int_0^L\psi_1\cos\frac{k\pi x}{L}dx\\
&-\xi \bar v (\frac{k\pi}{L})^2 \int_0^L\gamma_1\cos\frac{k\pi x}{L}dx=0,
\end{align}
and
\begin{align}\label{324}
\int_0^L\varphi_1\cos\frac{k\pi x}{L}dx+\int_0^L\psi_1\cos\frac{k\pi x}{L}dx- \Big(\Big(\frac{k\pi}{L}\Big)^2+\lambda \Big)\int_0^L\gamma_1\cos\frac{k\pi x}{L}dx=0.
\end{align}
On the other hand, since $(\varphi_1, \psi_1, \gamma_1)\in\mathcal{Z}$, we have from (\ref{316}) that
\begin{equation}\label{325}
P_k\int_0^L\varphi_1\cos\frac{k\pi x}{L}dx+Q_k\int_0^L\psi_1\cos\frac{k\pi x}{L}dx+\int_0^L\gamma_1\cos\frac{k\pi x}{L}dx=0,
\end{equation}
where $P_k$, $Q_k$ are defined in (\ref{311}) and (\ref{312}).  From (\ref{323})--(\ref{325}), we arrive at the following system
\begin{equation*}
\left(
  \begin{array}{ccc}
    - a_2\mu_2 \bar v&-(d_2(\frac{k\pi}{L})^2+\mu_2\bar v)&\xi\bar v (\frac{k\pi}{L})^2\\
1&1&-\big((\frac{k\pi}{L})^2+\lambda \big)\\
P_k&Q_k&1
\end{array}
\right)
\left(
  \begin{array}{ccc}
    \int_0^L\varphi_1\cos\frac{k\pi x}{L}dx\\ \int_0^L\psi_1\cos\frac{k\pi x}{L}dx\\
     \int_0^L\gamma_1\cos\frac{k\pi x}{L}dx
\end{array}
\right)=
\left(
  \begin{array}{ccc}
    0\\
0\\
0
\end{array}
\right).
\end{equation*}
We claim that the coefficient matrix $\mathcal{M}$ is nonsingular.  Indeed, we have from straightforward calculations that its determinant equals
\begin{align}\label{326}
\text{det}(\mathcal{M})=&d_2(\frac{k\pi}{L})^2+(1-a_2)\mu_2\bar v+P_k \Big(\Big(d_2\Big(\frac{k\pi}{L}\Big)^2+\mu_2\bar v \Big) \Big(\Big(\frac{k\pi}{L}\Big)^2+\lambda \Big) \nonumber\\
&-\xi\bar v\Big(\frac{k\pi}{L}\Big)^2 \Big)+Q_k \Big(\xi\bar v\Big(\frac{k\pi}{L}\Big)^2-a_2 \mu_2 \bar v \Big(\Big(\frac{k\pi}{L}\Big)^2+\lambda \Big)\Big)\nonumber\\
\hspace{2.2cm}=&d_2\Big(\frac{k\pi}{L}\Big)^2+(1-a_2)\mu_2 \bar v+\frac{ \big( \big(d_2(\frac{k\pi}{L})^2+\mu_2\bar v \big) \big((\frac{k\pi}{L}\big)^2+\lambda \big)-\xi\bar v(\frac{k\pi}{L})^2 \big)^2}{d_2(\frac{k\pi}{L})^2+(1-a_2)\mu_2\bar v} \nonumber\\
&+\frac{ \big(\xi\bar v(\frac{k\pi}{L})^2-a_2\mu_2 \bar v \big((\frac{k\pi }{L}\big)^2+\lambda \big) \big)^2}{d_2(\frac{k\pi}{L})^2+(1-a_2)\mu_2 \bar v}>0
\end{align}
where the last identity follows by direct substitutions of (\ref{311}) and (\ref{312}).  Hence we can quickly have that
\begin{equation}\label{327}
\int_0^L\varphi_1\cos\frac{k\pi x}{L}dx=\int_0^L\psi_1\cos\frac{k\pi x}{L}dx=\int_0^L\gamma_1\cos\frac{k\pi x}{L}dx=0,
\end{equation}
and $\mathcal{K}_1=0$ follows from (\ref{322}) and (\ref{327}).

\begin{figure}[h!]
        \centering
        \begin{subfigure}[b]{0.45\textwidth}
                \includegraphics[width=\textwidth]{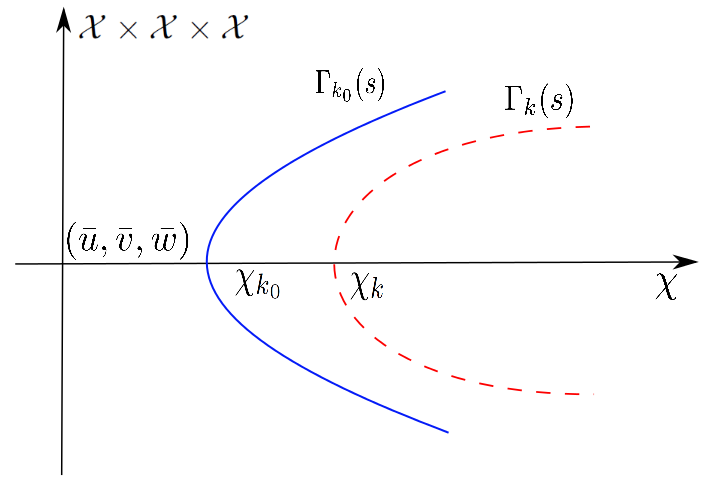}
                \caption*{If $\Gamma_{k_0}(s)$ turns to the right, the portion around $(\bar u,\bar v,\bar w,\chi_{k_0})$ is stable;  $\Gamma_{k}(s)$ around $(\bar u,\bar v,\bar w,\chi_{k})$ is always unstable if $k\neq k_0$.}
                \label{fig:gull}
        \end{subfigure}\hspace{0.15in}
        \begin{subfigure}[b]{0.45\textwidth}
        \includegraphics[width=\textwidth]{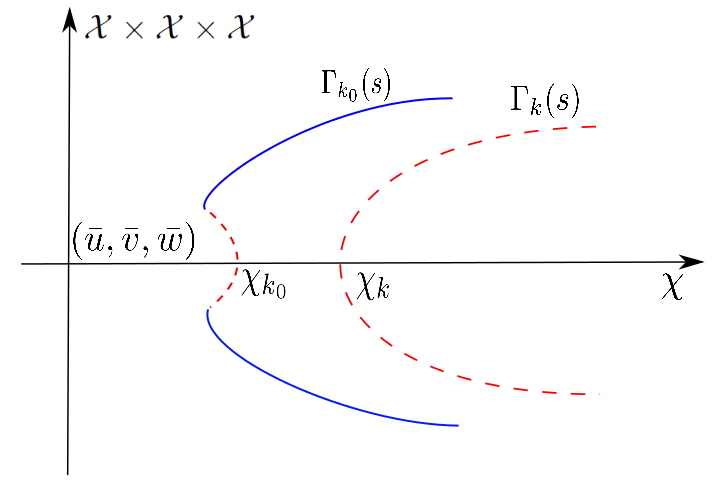}
                \caption*{If $\Gamma_{k_0}(s)$ turns to the left, the portion around $(\bar u,\bar v,\bar w,\chi_{k_0})$ is unstable;  $\Gamma_{k}(s)$ around $(\bar u,\bar v,\bar w,\chi_{k})$ is always unstable if $k\neq k_0$.}
                \label{fig:tiger}
        \end{subfigure}
 \caption{Bifurcation diagrams when case (i) of Proposition \ref{proposition2} occurs.  The stable bifurcation curve is plotted in solid line and the unstable bifurcation curve is plotted in dashed line.  The bifurcation curve is of pitch--fork type, i.e., being one--sided.}\label{fig1}
\end{figure}

Now we are ready to present our stability results which indicate that if $\mathcal{K}_2\neq 0$ then the sign of $\mathcal{K}_2$ determines the stability of the bifurcating solution in some cases.  To be precise, we have the following proposition.
\begin{proposition}\label{proposition2}
Suppose that all conditions in Theorem \ref{theorem31} are satisfied.  Let $(u_k(s,k),v_k(s,k), w_k(s,k), \chi_k(s))$ given by (\ref{315}) be a steady state of (\ref{31}), then $\mathcal{K}_1=0$ in (\ref{320}), i.e., the bifurcation is pitch--fork around $(\bar u,\bar v,\bar w,\chi_k)$, $\forall k\in \mathbb N^+$.  Denote $\chi_0=\min_{k\in \mathbb{N}^+}\{\chi_k,\hat \chi_k\}$ as in (\ref{33}).  Then we have the following stability results:

(i) If $\chi_0=\chi_{k_0}<\min_{k\in\mathbb N^+}\hat \chi_k$, then for all $k\neq k_0$, $(u_k(s,k),v_k(s,k), w_k(s,k), \chi_k(s))$ around $(\bar u,\bar v,\bar w,\chi_k)$ is always unstable, and $(u_{k_0}(s,k),v_{k_0}(s,k), w_{k_0}(s,k), \chi_{k_0}(s))$ around $(\bar u,\bar v,\bar w,\chi_{k_0})$ is asymptotically stable if $\mathcal{K}_2>0$ and it is unstable if $\mathcal{K}_2<0$.

(ii) If $\chi_0=\hat \chi_{k_1}<\min_{k\in\mathbb N^+} \chi_k$, then for all $k\in\mathbb N^+$, the bifurcating solutions $(u_k(s,k),v_k(s,k), w_k(s,k), \chi_k(s))$ around $(\bar u,\bar v,\bar w,\chi_k)$ are unstable.
\end{proposition}
The bifurcation curves $\Gamma_k(s)$ in case (i) are schematically illustrated in Figure \ref{fig1}.  Our results suggest that $(\bar u,\bar v,\bar w,\chi_k)$ loses its stability to stable steady state bifurcating solution with wavemode number $k_0$ where $\chi_k$ achieves its minimum over $\mathbb N^+$.  When case (ii) occurs, we surmise that the stability of the homogeneous solution is lost to stable Hopf bifurcating solutions.
\begin{remark}\label{remark3}
In general, it is not obvious to determine when case (i) or case (ii) occurs.  However when the interval length $L$ is sufficiently small, we have that
\[\tilde \chi_k \approx \frac{d_1}{\bar {u}}(\frac{k\pi}{L})^2< \frac{(d_1+d_2)^2+(d_1+d_2)d_1d_2+d_1+d_2}{(d_1+1)\bar {u}}(\frac{k\pi}{L})^2\approx \hat \chi_k\]
for all $k\in \mathbb N^+$.  This implies that for small intervals, $\chi_0=\min_{k\in\mathbb N^+} \chi_k=\chi_1$ and $(\bar u,\bar v,\bar w)$ loses its stability to the steady state bifurcating solution as $\chi$ surpasses $\chi_0$.  Since $k_0=1$, the bifurcating solution has a stable wavemode $\cos \frac{\pi x}{L}$ which is monotone in $x$.  Moreover, it is also easy to see that $k_0$ is increasing or non-decreasing in $L$.  Therefore, small domain only supports monotone stable solutions, while large interval supports non-monotone solutions, at least when $\chi$ is around $\chi_0$.  Actually, if $(u(x),v(x),w(x))$ is an increasing solution to (\ref{31}), $(u(L-x),v(L-x),w(L-x))$ is a decreasing solution, then one can construct non-monotone solutions to (\ref{31}) over $(0,2L)$, $(0,3L)$,...by reflecting and periodically extending the monotone ones at $x=L$, $2L$, $3L$,...
\end{remark}
\begin{proof} [Proof\nopunct] \emph{of Proposition} \ref{proposition2}.
We shall only prove case (ii), while the unstable part in case (i) can be proved by the same arguments and the stable part in case (i) can be proved following slight modification in the arguments for Corollary 1.13 of \cite{CR2}, or Theorem 3.2 of \cite{WGY}, Theorem 5.5, Theorem 5.6 of \cite{CKWW}.

For each $k\in \mathbb N^+$, we linearize (\ref{31}) around $(u_k(s,x),v_k(s,x),w_k(s,x),\chi_k(s))$ for $s\in (-\delta,\delta)$ and obtain the following eigenvalue problem
\[D_{(u,v,w)}\mathcal{F}(u_k(s,x),v_k(s,x),w_k(s,x),\chi_k(s))(u,v,w)={\sigma_k}(s)(u,v,w),~(u,v,w)\in \mathcal{X} \times \mathcal{X}\times \mathcal{X}.\]
$(u_k(s,x),v_k(s,x),w_k(s,x),\chi_k(s))$ is asymptotically stable if and only if the real part of any eigenvalue ${\sigma_k}(s)$ is negative.

Sending $s\rightarrow 0$, we know from bifurcation analysis that ${\sigma_k}=0$ is a simple eigenvalue of $D_{(u,v,w)}\mathcal{F}(\bar{u},\bar{v},\bar{w},\chi_k)={\sigma_k}(u,v,w)$ or equivalently, the following eigenvalue problem
\begin{eqnarray*}
\left\{
\begin{array}{ll}
d_1u''-\chi_k\bar uw''-\mu_1\bar uu-\mu_1a_1\bar uv={\sigma_k} u,&x\in(0,L),\\
d_2v''-\xi\bar vw''-\mu_2\bar vv-\mu_2a_2\bar vu={\sigma_k} v,&x\in(0,L), \\
w''-\lambda w+u+v={\sigma_k} w,&x\in(0,L),\\
u'(x)=v'(x)=w'(x)=0,&x=0,L;
\end{array}
\right.
\end{eqnarray*}
moreover, it has a one--dimensional eigen-space $\mathcal{N}\big(D_{(u,v,w)}\mathcal{F}(\bar{u},\bar{v},\bar{w},\chi_k)\big)=\{(P_k,Q_k,1)\cos \frac{k\pi x}{L}\}$ and it is easy to show that $(P_k,Q_k,1)\cos \frac{k\pi x}{L} \not\in\mathcal{R}\big(D_{(u,v,w)}\mathcal{F}(\bar{u},\bar{v},\bar{w},\chi_k)\big)$ following the same analysis that leads to (\ref{317}).  Multiplying the system  above by $\cos \frac{k\pi x}{L}$ and integrating them over $(0,L)$ by parts, we have that ${\sigma_k}$ is an eigenvalue of (\ref{32}) with $\chi=\chi_k$ or the following matrix
\[
\begin{pmatrix}
-(d_1(\frac{k\pi}{L})^2+\mu_1\bar u) & -a_1\mu_1\bar u & \chi_k\bar u(\frac{k\pi}{L})^2 \\
-a_2\mu_2\bar v & -(d_2(\frac{k\pi}{L})^2+\mu_2\bar v) & \xi\bar v(\frac{k\pi}{L})^2\\
1& 1& -((\frac{k\pi}{L})^2+\lambda)
\end{pmatrix}.
\]
If $\chi_0=\min_{k\in\mathbb N^+}\hat \chi_k< \chi_k$ for all $k\in\mathbb N^+$, or $\chi_0=\min_{k\in\mathbb N^+}\chi_k< \hat \chi_k$ for $k \neq k_0$, we have from the proof of Proposition \ref{proposition1} that the matrix above always has an eigenvalue ${\sigma_k}$ with positive real part in each case.  From the standard eigenvalue perturbation theory in \cite{Ka}, for $s$ being small, there exists an eigenvalue ${\sigma_k}(s)$ to the linearized problem above that has a positive real part.  This proves the instability of $(u_k(s,x),v_k(s,x),w_k(s,x),\chi_k(s))$ for $s\in(-\delta,\delta)$.
\end{proof}

We now proceed to evaluate $\mathcal{K}_2$ which determines the turning direction of $\Gamma_k(s)$ and the stability of $\Gamma_{k_0}(s)$ for $s\in(-\delta,\delta)$.  Equating the $s^3$--terms in (\ref{31}), we collect the following system
\begin{equation}\label{328}
\left\{
\begin{array}{ll}
d_1\varphi_{2}^{''}\!-\!\chi_k \bar u\gamma_{2}^{''}\!=\!\chi_k \Big(\! -(\frac{k\pi}{L}\varphi_{1}^{'}\!+\!P_k\frac{k\pi}{L}\gamma_{1}^{'})\sin\frac{k\pi x}{L}\!+\!\big(\!-(\frac{k\pi }{L})^2\varphi_1\!+\!P_k\gamma_{1}^{''}\!\big)\cos\frac{k\pi x}{L}\!\Big)\\
\hspace{2.7cm}+\mu_1 \Big(\big((2P_k+a_1 Q_k)\varphi_1+a_1P_k\psi_{1}\big)\cos\frac{k\pi x}{L}+\bar u(\varphi_2+a_1\psi_2) \Big)\\
\hspace{2.7cm}-\bar u (\frac{k\pi}{L})^2\mathcal{K}_2\cos\frac{k\pi x}{L},\\
d_2\psi_{2}^{''}-\xi \bar v\gamma_{2}^{''}=\xi \Big(-(\frac{k\pi}{L}\psi_{1}^{'}+Q_k\frac{k\pi}{L}\gamma_{1}^{'})\sin\frac{k\pi x}{L}+\big(-(\frac{k\pi}{L})^2\psi_1+Q_k\gamma_{1}^{''} \big)\cos\frac{k\pi x}{L} \Big)\\
\hspace{2.5cm}+\mu_2 \Big((a_2Q_k\varphi_{1}+(a_2P_k+2Q_k)\psi_1)\cos\frac{k\pi x}{L}+\bar v(a_2\varphi_2+\psi_2) \Big),\\
\gamma_{2}^{''}-\lambda\gamma_2+\varphi_2+\psi_2 =0,\\
\varphi_{2}^{'}(x)=\psi_{2}^{'}(x)=\gamma_{2}^{'}(x)=0, x=0, L.
\end{array}
\right.
\end{equation}
Multiplying the first equation in (\ref{328}) by $\cos\frac{k\pi x}{L}$, we conclude from the integration by parts and straightforward calculations that
\begin{align}\label{329}
\frac{\bar u(k\pi)^2}{2L}\mathcal{K}_2=&\Big(\Big(\frac{d_1(\frac{k\pi}{L})^2}{2\bar u}+\frac{3\mu_1}{2}\Big)P_k +a_1\mu_1Q_k \Big)\int_0^L\varphi_1\cos\frac{2k\pi x}{L}dx \nonumber\\
&+\frac{a_1\mu_1}{2}P_k\int_0^L\psi_1\cos\frac{2k\pi x}{L}dx-\Big(\Big(\frac{d_1(\frac{k\pi}{L})^2}{\bar u}+\mu_1\Big)P_{k}^{2} +a_1\mu_1P_kQ_k \Big)\nonumber\\
&\cdot\int_0^L\gamma_1\cos\frac{2k\pi x}{L}dx+\Big(d_1\Big(\frac{k\pi}{L}\Big)^2+\mu_1\bar u \Big)\int_0^L\varphi_2\cos\frac{k\pi x}{L}dx\\
&+a_1\mu_1\bar u\int_0^L\psi_2\cos\frac{k\pi x}{L}dx-\Big(\Big(d_1\Big(\frac{k\pi}{L}\Big)^2+\mu_1\bar u\Big)P_k+a_1\mu_1\bar u Q_k\Big)\nonumber\\
&\cdot\int_0^L\gamma_2\cos\frac{k\pi x}{L}dx-\Big(\frac{d_1(\frac{k\pi}{L})^2}{2\bar u}-\frac{\mu_1}{2}\Big)P_k\int_0^L\varphi_1dx+\frac{a_1\mu_1}{2}P_k\int_0^L\psi_1dx, \nonumber
\end{align}
where we have substituted $\chi_k$ given by (\ref{39}) and used the fact that $P_k+Q_k=(\frac{k\pi}{L})^2+\lambda$.  On the other hand, we can test the second and the third equation of (\ref{328}) by $\cos\frac{k\pi x}{L}$ over $(0,L)$ to obtain
\begin{align}\label{330}
&-a_2\mu_2\bar v\int_0^L\varphi_2\cos\frac{k\pi x}{L}dx-\Big(d_2\Big(\frac{k\pi}{L}\Big)^2+\mu_2\bar v \big)\int_0^L\psi_2\cos\frac{k\pi x}{L}dx \nonumber\\
&+\xi\bar v\Big(\frac{k\pi}{L}\Big)^2\int_0^L\gamma_2\cos\frac{k\pi x}{L}dx \nonumber\\
=&\frac{a_2\mu_2}{2}Q_k\int_0^L\varphi_1\cos\frac{2k\pi x}{L}dx+\Big(\frac{\xi}{2}\Big(\frac{k\pi}{L}\Big)^2+\frac{a_2\mu_2}{2}P_k+\mu_2Q_k \Big)\int_0^L\psi_1\cos\frac{2k\pi x}{L}dx \nonumber\\
&-\xi Q_k\Big(\frac{k\pi}{L}\Big)^2\int_0^L\gamma_1\cos\frac{2k\pi x}{L}dx+\frac{a_2\mu_2}{2}Q_k\int_0^L\varphi_1dx \nonumber\\
&+\Big(\frac{a_2\mu_2}{2}P_k+\mu_2Q_k-\frac{\xi}{2}\Big(\frac{k\pi}{L}\Big)^2 \Big)\int_0^L\psi_1dx,
\end{align}
and
\begin{align}\label{331}
\int_0^L\varphi_2\cos\frac{k\pi x}{L}dx+\int_0^L\psi_2\cos\frac{k\pi x}{L}dx-\Big(\Big(\frac{k\pi}{L}\Big)^2+\lambda \Big)\int_0^L\gamma_2\cos\frac{k\pi x}{L}dx=0.
\end{align}
Moreover, it follows from $(\varphi_2, \psi_2, \gamma_2)\in \mathcal{Z}$ defined in (\ref{316}) that,
\begin{align}\label{332}
P_k\int_0^L\varphi_2\cos\frac{k\pi x}{L}dx+Q_k\int_0^L\psi_2\cos\frac{k\pi x}{L}dx+\int_0^L\gamma_2\cos\frac{k\pi x}{L}dx=0.
\end{align}
Putting (\ref{330})--(\ref{332}) together, we arrive at the following system
\begin{equation}\label{333}
\left(
  \begin{array}{ccc}
    \!\!\!-a_2\mu_2\bar v&-(\!d_2(\!\frac{k\pi}{L}\!)^2\!+\!\mu_2\bar v\!)&\!\xi\bar v (\!\frac{k\pi}{L}\!)^2\\
\!\!\!1&\!1&\!-\big(\!(\!\frac{k\pi}{L}\!)^2\!+\!\lambda \!\big)\\
\!\!\!P_k&\!Q_k&\!1
\end{array}
\!\!\right)
\!\!\left(
  \begin{array}{ccc}
    \!\!\!\int_0^L\!\varphi_2\cos\frac{k\pi x}{L}dx\\  \!\int_0^L\!\psi_2\cos\frac{k\pi x}{L}dx\\  \!\int_0^L\!\gamma_2\cos\frac{k\pi x}{L}dx
\end{array}
\!\!\right)\!\!=\!\!
\left(\!\!
  \begin{array}{ccc}
   G\\
0\\
0
\end{array}
\!\!\right),
\end{equation}
where we have used the notation
\begin{align}\label{334}
G=&\frac{a_2\mu_2}{2}Q_k\int_0^L\varphi_1\cos\frac{2k\pi x}{L}dx+\Big(\frac{\xi}{2}\Big(\frac{k\pi}{L}\Big)^2+\frac{a_2\mu_2}{2}P_k+\mu_2Q_k \Big)\int_0^L\psi_1\cos\frac{2k\pi x}{L}dx \nonumber\\
&-\xi Q_k\Big(\frac{k\pi}{L}\Big)^2\int_0^L\gamma_1\cos\frac{2k\pi x}{L}dx+\frac{a_2\mu_2}{2}Q_k\int_0^L\varphi_1dx \\
&+\Big(\frac{a_2\mu_2}{2}P_k+\mu_2Q_k-\frac{\xi}{2}\Big(\frac{k\pi}{L}\Big)^2 \Big) \int_0^L\psi_1dx.\nonumber
\end{align}
Solving the algebraic system (\ref{333}) gives us
\begin{equation}\label{335}
\int_0^L\varphi_2\cos\frac{k\pi x}{L}dx=\frac{\vert \mathcal{A}_1\vert}{\vert \mathcal{A}\vert}~,~\int_0^L\psi_2\cos\frac{k\pi x}{L}dx=\frac{\vert \mathcal{A}_2\vert}{\vert \mathcal{A}\vert},
\end{equation}
and
\begin{equation}\label{336}
\int_0^L\gamma_2\cos\frac{k\pi x}{L}dx=\frac{\vert \mathcal{A}_3\vert}{\vert \mathcal{A}\vert},
\end{equation}
where we have denoted
\begin{equation}\label{337}
\mathcal{A}=\left(
  \begin{array}{ccc}
    -a_2\mu_2 \bar v&-\big(d_2(\frac{k\pi}{L})^2+\mu_2\bar v\big)&\xi\bar v (\frac{k\pi}{L})^2\\
1&1&-\big((\frac{k\pi}{L})^2+\lambda \big)\\
P_k&Q_k&1
\end{array}
\right),
\end{equation}
\begin{equation}\label{338}
\mathcal{A}_1=\left(
  \begin{array}{ccc}
    G&-\big(d_2(\frac{k\pi}{L})^2+\mu_2\bar v\big)&\xi\bar v (\frac{k\pi}{L})^2\\
0&1&-\big((\frac{k\pi}{L})^2+\lambda \big)\\
0&Q_k&1
\end{array}
\right),
\end{equation}
\begin{equation}\label{339}
\mathcal{A}_2=\left(
  \begin{array}{ccc}
    -a_2\mu_2 \bar v&G&\xi\bar v (\frac{k\pi}{L})^2\\
1&0&-\big((\frac{k\pi}{L})^2+\lambda \big)\\
P_k&0&1
\end{array}
\right),
\end{equation}
and
\begin{equation}\label{340}
\mathcal{A}_3=\left(
  \begin{array}{ccc}
  -a_2\mu_2\bar v&-\big(d_2(\frac{k\pi}{L})^2+\mu_2\bar v\big)&G\\
1&1&0\\
P_k&Q_k&0
\end{array}
\right).
\end{equation}

Integrating (\ref{321}) over $(0,L)$ by parts, thanks to the non-flux boundary conditions and $\mathcal{K}_1=0$, we obtain that
\begin{equation*}
\left\{
\begin{array}{ll}
\int_0^L\varphi_1dx+a_1\int_0^L\psi_1dx=-\frac{L}{2\bar u}P_k(P_k+a_1Q_k) \\
a_2\int_0^L\varphi_1dx+\int_0^L\psi_1dx=-\frac{L}{2\bar v}Q_k(Q_k+a_2P_k),\\
\int_0^L\varphi_1dx+\int_0^L\psi_1dx-\lambda\int_0^L\gamma_1dx=0,
\end{array}
\right.
\end{equation*}
 or equivalently
\begin{equation}\label{341}
\left(
  \begin{array}{ccc}
   1&a_1&0\\
a_2&1&0\\
1&1&-\lambda
\end{array}
\right)
\left(
  \begin{array}{ccc}
    \int_0^L\varphi_1dx\\ \int_0^L\psi_1dx\\ \int_0^L\gamma_1dx
\end{array}
\right)=
\left(
  \begin{array}{ccc}
  -\frac{L}{2\bar u}P_k(P_k+a_1Q_k)\\
  -\frac{L}{2\bar v}Q_k(Q_k+a_2P_k)\\
0
\end{array}
\right).
\end{equation}
Solving (\ref{341}) leads us to
\begin{equation}\label{342}
\int_0^L\varphi_1dx=\frac{\vert \mathcal{B}_1\vert}{\vert \mathcal{B}\vert},~\int_0^L\psi_1dx=\frac{\vert \mathcal{B}_2\vert}{\vert \mathcal{B}\vert},~\int_0^L\gamma_1 dx=\frac{\vert \mathcal{B}_3\vert}{\vert \mathcal{B}\vert},
\end{equation}
where
\begin{equation}\label{343}
\mathcal{B}=\left(
  \begin{array}{ccc}
   1&a_1&0\\
a_2&1&0\\
1&1&-\lambda
\end{array}
\right),
\end{equation}
\begin{equation}\label{344}
\mathcal{B}_1=\left(
  \begin{array}{ccc}
    -\frac{L}{2\bar u}P_k(P_k+a_1Q_k)&a_1&0\\
-\frac{L}{2\bar v}Q_k(Q_k+a_2P_k)&1&0\\
0&1&-\lambda
\end{array}
\right),
\end{equation}
\begin{equation}\label{345}
\mathcal{B}_2=\left(
  \begin{array}{ccc}
 1&-\frac{L}{2\bar u}P_k(P_k+a_1Q_k)&0\\
a_2&-\frac{L}{2\bar v}Q_k(Q_k+a_2P_k)&0\\
1&0&-\lambda
\end{array}
\right),
\end{equation}
and
\begin{equation}\label{346}
\mathcal{B}_3=\left(
  \begin{array}{ccc}
 1&a_1&-\frac{L}{2\bar u}P_k(P_k+a_1Q_k)\\
a_2&1&-\frac{L}{2\bar v}Q_k(Q_k+a_2P_k)\\
1&1&0
\end{array}
\right).
\end{equation}

Similarly, we test (\ref{321}) by $\cos\frac{2k\pi x}{L}$ and obtain that
\begin{align}\label{347}
&\left(
  \begin{array}{ccc}
   \!\!\!-\big(d_1(\frac{2k\pi}{L})^2+\mu_1 \bar u \big)\!&-a_1\mu_1\bar u\!&\chi_k\bar u(\frac{2k\pi}{L})^2\\
\!\!\!-a_2\mu_2 \bar v\!&-\big(d_2(\frac{2k\pi}{L})^2+\mu_2 \bar v\big)\!&\xi\bar v(\frac{2k\pi}{L})^2\\
\!\!\!1\!&1\!&-\big((\frac{2k\pi}{L})^2+\lambda \big)
\end{array}
\!\!\!\right)
\left(
  \begin{array}{ccc}
    \!\!\!\int_0^L\varphi_1\cos\frac{2k\pi x}{L}dx\\ \!\!\!\int_0^L\psi_1\cos\frac{2k\pi x}{L}dx\\ \!\!\!\int_0^L\gamma_1\cos\frac{2k\pi x}{L}dx
\end{array}
\!\!\!\right) \nonumber\\
=&\left(
  \begin{array}{ccc}
 -\frac{\chi_k L}{2}P_k(\frac{k\pi}{L})^2+\frac{\mu_1 L}{4}P_k(P_k+a_1Q_k)\\
-\frac{\xi L}{2}Q_k(\frac{k\pi}{L})^2+\frac{\mu_2 L}{4}Q_k(Q_k+a_2P_k)\\
0
\end{array}
\right).
\end{align}
Solving system (\ref{347}), we have that
\begin{equation}\label{348}
\int_0^L\varphi_1\cos\frac{2k\pi x}{L}dx=\frac{\vert \mathcal{C}_1\vert}{\vert \mathcal{C}\vert}~,~\int_0^L\psi_1\cos\frac{2k\pi x}{L}dx=\frac{\vert \mathcal{C}_2\vert}{\vert \mathcal{C}\vert},
\end{equation}
and
\begin{equation}\label{349}
\int_0^L\gamma_1\cos\frac{2k\pi x}{L} dx=\frac{\vert \mathcal{C}_3\vert}{\vert \mathcal{C}\vert},
\end{equation}
where we have applied the following notations in (\ref{348}) and (\ref{349}),
\begin{equation}\label{350}
\mathcal{C}=\left(
  \begin{array}{ccc}
  -\big(d_1(\frac{2k\pi}{L})^2+\mu_1 \bar u\big)&-a_1\mu_1\bar u&\chi_k\bar u(\frac{2k\pi}{L})^2\\
-a_2\mu_2 \bar v&-\big(d_2(\frac{2k\pi}{L})^2+\mu_2 \bar v\big)&\xi\bar v(\frac{2k\pi}{L})^2\\
1&1&-\big((\frac{2k\pi}{L})^2+\lambda \big)
\end{array}
\right),
\end{equation}
\begin{equation}\label{351}
\mathcal{C}_1\!\!=\!\!\left(
  \begin{array}{ccc}
    \!\!\!\!-\frac{\chi_k L}{2} P_k(\!\frac{k\pi}{L}\!)^2\!+\!\frac{\mu_1 L}{4}P_k (\!P_k\!+\!a_1Q_k\!)\!\!&-a_1\mu_1  \bar u\!\!&\chi_k\bar u (\frac{2k\pi}{L})^2\\
\!\!\!\!-\frac{\xi L}{2}Q_k(\!\frac{k\pi}{L}\!)^2\!+\!\frac{\mu_2 L}{4}Q_k (\!Q_k\!+\!a_2P_k\!)\!\!&-\big(\!d_2(\!\frac{2k\pi}{L}\!)^2\!+\!\mu_2 \bar v\!\big)\!\!&\xi\bar v(\!\frac{2k\pi}{L}\!)^2\\
\!\!\!\!0\!\!&1\!\!&-\big(\!(\!\frac{2k\pi}{L}\!)^2+\lambda\! \big)
\end{array}
\!\!\!\!\right),
\end{equation}
\begin{equation}\label{352}
\mathcal{C}_2\!\!=\!\!\left(
  \begin{array}{ccc}
 \!\!\!\!-\big(\!d_1(\!\frac{2k\pi}{L}\!)^2\!+\!\mu_1 \bar u\!\big)& \!-\frac{\chi_k L}{2} P_k(\!\frac{k\pi}{L}\!)^2\!+\!\frac{\mu_1 L}{4}P_k (\!P_k\!+\!a_1Q_k\!)&\!\chi_k\bar u (\!\frac{2k\pi}{L}\!)^2\\
\!\!\!\!-a_2\mu_2 \bar v&\!-\frac{\xi L}{2}Q_k(\!\frac{k\pi}{L}\!)^2\!+\!\frac{\mu_2 L}{4}Q_k(\!Q_k\!+\!a_2P_k\!)&\!\xi\bar v(\!\frac{2k\pi}{L}\!)^2\\
\!\!\!\!1&\!0&\!-\big(\!(\!\frac{2k\pi}{L}\!)^2+\lambda \!\big)
\end{array}
\!\!\!\!\right),
\end{equation}
and
\begin{equation}\label{353}
\!\!\!\!\!\mathcal{C}_3\!\!=\!\!\left(
  \begin{array}{ccc}
 \!\!\!\!-\big(\!d_1(\!\frac{2k\pi}{L}\!)^2\!+\!\mu_1 \bar u\!\big)&\!\!\!\!\!-a_1\mu_1 \bar u&\!\!\!\!\!-\frac{\chi_k L}{2}P_k (\!\frac{k\pi}{L}\!)^2\!+\!\frac{\mu_1 L}{4}P_k(\!P_k\!+\!a_1Q_k\!)\\
\!\!\!\!-a_2\mu_2 \bar v&\!\!\!\!\!-\big(\!d_2(\!\frac{2k\pi}{L})^2\!+\!\mu_2 \bar v\!\big)&\!\!\!\!\!-\frac{\xi L}{2}Q_k(\!\frac{k\pi}{L}\!)^2\!+\!\frac{\mu_2 L}{4}Q_k(\!Q_k\!+\!a_2P_k\!)\\
\!\!\!\!1&\!\!\!\!\!1&\!\!\!\!\!0
\end{array}
\!\!\!\!\right).
\end{equation}%

We recall from (\ref{329}) that $\mathcal{K}_2$ consists of integrals (\ref{335})--(\ref{336}), (\ref{342}) and (\ref{348})--(\ref{349}).  As we can see from the formulas above, the calculations to evaluate $\mathcal{K}_2$ are extremely complicated and lengthy, therefore, for the simplicity of calculations, we suppose that both $d_1$ and $d_2$ are large throughout the rest of our calculations.  Therefore, we can perform precise calculations to evaluate the sign of $\mathcal K_2$ hence the stability of the bifurcating solutions with $k=k_0$ in case (i). Biologically, this interprets the scenario that both species $u$ and $v$ diffuse much faster than the chemical.  According to the results of Tello and Winkler \cite{TW}, when the chemotaxis coefficients $\chi$ and $\xi$ are small, $(\bar u,\bar v,\bar w)$ is a global attractor of (\ref{11}) and it prevents the existence of nonconstant positive solutions to (\ref{31}).  On the other hand, it is well known that the dynamics of (\ref{11}) is dominated by that of the ODE if either $d_1$ or $d_2$ is sufficiently large.  However, when the diffusion rates are large, we can always choose even larger chemotaxis rate $\chi$ to support the formation of nontrivial and interesting patterns.  This is numerically verified in Section \ref{section4}.

First of all, we have the following asymptotic expansions from (\ref{39}) and (\ref{311})--(\ref{312}) as $d_1\rightarrow \infty$ with $d_2=O(d_1)$
\begin{equation}\label{354}
\chi_k=\Big(\frac{(\frac{k\pi}{L})^4((\frac{k\pi}{L})^2+\lambda)}{(\frac{k\pi}{L})^4\bar u}+O(\frac{1}{d_1})\Big)d_1=\frac{(\frac{k\pi}{L})^2+\lambda}{\bar u}d_1+O(1),
\end{equation}
\begin{equation}\label{355}
P_k=(\frac{k\pi}{L})^2+\lambda+O(\frac{1}{d_2}),Q_k=O(\frac{1}{d_2}).
\end{equation}
These asymptotic formulas will be used frequently throughout the rest of this section.

Substituting (\ref{354})--(\ref{355}) into (\ref{350}) and (\ref{351}), we collect
\begin{align}\label{356}
\vert\mathcal{C}\vert=&-\big(d_1(\frac{2k\pi}{L})^2+\mu_1\bar u \big)\big(d_2(\frac{2k\pi}{L})^2+\mu_2\bar v \big)\big((\frac{2k\pi}{L})^2+\lambda \big)-\chi_k a_2\mu_2 \bar u\bar v(\frac{2k\pi}{L})^2\nonumber\\
&+\chi_k\bar u(\frac{2k\pi}{L})^2 \big(d_2(\frac{2k\pi}{L})^2+\mu_2\bar v \big)+a_1a_2\mu_1\mu_2\bar u\bar v \big((\frac{2k\pi}{L})^2+\lambda \big)\nonumber\\
&+\xi\bar v(\frac{2k\pi}{L})^2 \big(d_1(\frac{2k\pi}{L})^2+\mu_1\bar u \big)-\xi a_1\mu_1\bar u\bar v(\frac{2k\pi}{L})^2\nonumber\\
=&\Big(-(\frac{2k\pi}{L})^4\Big((\frac{2k\pi}{L})^2+\lambda\Big)+\bar u(\frac{2k\pi}{L})^4\frac{(\frac{k\pi }{L})^2+\lambda}{\bar u}+o(1)\Big)d_1d_2\nonumber\\
=&\big(-48(\frac{k\pi}{L})^6+o(1)\big)d_1d_2,
\end{align}
and
\begin{align}\label{357}
\!\vert\mathcal{C}_1\vert\!\!=&\!\Big( \!\!\big(\!d_2(\!\frac{2k\pi}{L}\!)^2\!+\!\mu_2\bar v \!\big)\big(\!(\!\frac{2k\pi}{L}\!)^2\!+\!\lambda \!\big)\!-\!\xi\bar v(\!\frac{2k\pi}{L}\!)^2 \!\Big) \big(\!-\frac{\chi_kP_k L}{2}(\!\frac{k\pi}{L}\!)^2\!+\!\frac{\mu_1 L}{4}P_k(\!P_k\!+\!a_1Q_k\!) \!\big)\nonumber\\
&+\Big( \!\chi_k\bar u(\frac{2k\pi}{L})^2\!-\!a_1\mu_1\bar u \big((\frac{2k\pi}{L})^2\!+\!\lambda \big) \!\Big) \big(-\frac{\xi Q_k L}{2}(\frac{k\pi}{L})^2\!+\!\frac{\mu_2 L}{4}Q_k(Q_k\!+\!a_2P_k) \big)\nonumber\\
=&\Big(-2L(\frac{k\pi}{L})^4((\frac{2k\pi}{L})^2+\lambda)\frac{\chi_k P_k}{d_1}+o(1)\Big)d_1d_2\nonumber\nonumber\\
=&\Big(-2L(\frac{k\pi}{L})^4\Big((\frac{2k\pi}{L})^2+\lambda\Big)\frac{((\frac{k\pi}{L})^2+\lambda)^2}{\bar u}+o(1)\Big)d_1d_2.
\end{align}
Together with (\ref{356})--(\ref{357}), we have from straightforward calculations that
\begin{align}\label{358}
&\Big(\Big(\frac{d_1(\frac{k\pi}{L})^2}{2\bar u}+\frac{3\mu_1}{2}\Big)P_k +a_1\mu_1Q_k \Big)\int_0^L\varphi_1\cos\frac{2k\pi x}{L}dx \nonumber\\
=&\Big(\frac{(d_1(\frac{k\pi}{L})^2+3\mu_1\bar u)((\frac{k\pi}{L})^2+\lambda)}{2\bar u}+O(\frac{1}{d_1})\Big)\frac{\vert\mathcal{C}_1\vert}{\vert\mathcal{C}\vert}\nonumber\\
=&\Big(\frac{(d_1(\frac{k\pi}{L})^2+3\mu_1\bar u)((\frac{k\pi}{L})^2+\lambda)}{2\bar u}+O(\frac{1}{d_1})\Big)\cdot\frac{-2L(\frac{k\pi}{L})^4\frac{((\frac{2k\pi}{L})^2+\lambda)((\frac{k\pi}{L})^2+\lambda)^2}{\bar u}+o(1)}{-48(\frac{k\pi}{L})^6+o(1)}\nonumber\\
=&\frac{L((\frac{k\pi}{L})^2+\lambda)^3((\frac{2k\pi}{L})^2+\lambda)}{48{\bar u}^2}d_1+O(1).
\end{align}
We are ready to present the stability or instability of the nonconstant positive bifurcating steady states $(u_k(s,x),v_k(s,x),w_k(s,x))$ constructed in Theorem \ref{theorem31}.

\begin{theorem}\label{theorem32}
Suppose that all conditions in Theorem \ref{theorem31} are satisfied and denote $ \chi_0=\min_{k\in \mathbb{N}^+}\{\chi_k,\hat \chi_k\}$ as in (\ref{33}).  The bifurcating solutions $(u_k(s,x),v_k(s,x),w_k(s,x), \chi_k(s))$ around $(\bar u,\bar v,\bar w,\chi_k)$ is always unstable for each $k\in\mathbb N^+$ if $\chi_0=\hat \chi_{k_1}<\min_{k\in\mathbb N^+} \chi_k$ and is unstable for each positive integer $k\neq k_0$ if $\chi_0=\chi_{k_0}<\min_{k\in\mathbb N^+}\hat \chi_k$.  If $\chi_0=\chi_{k_0}$, $(u_{k_0}(s,k),v_{k_0}(s,k), w_{k_0}(s,k), \chi_{k_0}(s))$ around $(\bar u,\bar v,\bar w,\chi_{k_0})$ is asymptotically stable if $\mathcal{K}_2>0$ and it is unstable if $\mathcal{K}_2<0$.  In particular, there exists $\bar D$ large such that for all $d_1$ and $d_2=O(d_1)>\bar D$, $\mathcal{K}_2>0$ if $\lambda<\frac{(14-2a_1a_2)(\frac{k_0\pi}{L})^2}{1-a_1a_2}$ and $\mathcal{K}_2<0$ if $\lambda>\frac{(14-2a_1a_2)(\frac{k_0\pi}{L})^2}{1-a_1a_2}$.
\end{theorem}

\begin{proof}
Following the same calculations that lead to (\ref{358}), we obtain from (\ref{350}), $(\ref{352})$ and (\ref{353}) that
\begin{equation}\label{359}
\frac{a_1\mu_1}{2}P_{k_0}\int_0^L\psi_1\cos\frac{2{k_0}\pi x}{L}dx=\frac{a_1\mu_1P_{k_0}}{2}\frac{\vert\mathcal{C}_2\vert}{\vert\mathcal{C}\vert}=O(\frac{1}{d_2})=O(\frac{1}{d_1}),
\end{equation}
and
\begin{align}\label{360}
&\Big(\Big(\frac{d_1(\frac{{k_0}\pi}{L})^2}{\bar u}+\mu_1\Big)P_{{k_0}}^{2} +a_1\mu_1P_{k_0}Q_{k_0} \Big)\int_0^L\gamma_1\cos\frac{2{k_0}\pi x}{L}dx\nonumber\\
=&\Big(\Big(\frac{d_1(\frac{{k_0}\pi}{L})^2}{\bar u}+\mu_1\Big)P_{{k_0}}^{2} +a_1\mu_1P_{k_0}Q_{k_0} \Big)\frac{\vert\mathcal{C}_3\vert}{\vert\mathcal{C}\vert}=\frac{L((\frac{{k_0}\pi}{L})^2+\lambda)^4}{24{\bar u}^2}d_1+O(1).
\end{align}
Moreover, it follows from (\ref{337})--(\ref{340}) that
\begin{equation}\label{361}
\Big(d_1\Big(\frac{{k_0}\pi}{L}\Big)^2+\mu_1\bar u \Big)\int_0^L\varphi_2\cos\frac{{k_0}\pi x}{L}dx=\Big(d_1\Big(\frac{{k_0}\pi}{L}\Big)^2+\mu_1\bar u \Big)\frac{\vert\mathcal{A}_1\vert}{\vert\mathcal{A}\vert}=O(1),
\end{equation}
\begin{align}\label{362}
&a_1\mu_1\bar u\int_0^L\psi_2\cos\frac{{k_0}\pi x}{L}dx=a_1\mu_1\bar u\frac{\vert\mathcal{A}_2\vert}{\vert\mathcal{A}\vert}=o(\frac{1}{d_2}),
\end{align}
and
\begin{align}\label{363}
&\Big(\Big(d_1\Big(\frac{{k_0}\pi}{2}\Big)^2+\mu_1\bar u\Big)P_{k_0}+a_1\mu_1\bar u Q_{k_0}\Big)\int_0^L\gamma_2\cos\frac{{k_0}\pi x}{L}dx\nonumber\\
=&\Big(\Big(d_1\Big(\frac{{k_0}\pi}{2}\Big)^2+\mu_1\bar u\Big)P_{k_0}+a_1\mu_1\bar u Q_{k_0}\Big)\frac{\vert\mathcal{A}_3\vert}{\vert\mathcal{A}\vert}=O(1).
\end{align}

Finally, we conclude from (\ref{343})--(\ref{345}) that
\begin{align}\label{364}
&\Big(\frac{d_1(\frac{{k_0}\pi}{L})^2}{2\bar u}-\frac{\mu_1}{2}\Big)P_{k_0}\int_0^L\varphi_1dx=\Big(\frac{d_1(\frac{{k_0}\pi}{L})^2}{2\bar u}-\frac{\mu_1}{2}\Big)P_{k_0}\frac{\vert\mathcal{B}_1\vert}{\vert\mathcal{B}\vert}\nonumber\\
=&\frac{L((\frac{{k_0}\pi}{L})^2+\lambda)^3}{4(1-a_1a_2){\bar u}^2}(\frac{{k_0}\pi}{L})^2d_1+O(1),
\end{align}
and
\begin{align}\label{365}
&\frac{a_1\mu_1}{2}P_{k_0}\int_0^L\psi_1dx=\frac{a_1\mu_1}{2}P_{k_0}\frac{\vert\mathcal{B}_2\vert}{\vert\mathcal{B}\vert}=O(1).
\end{align}

Putting (\ref{358})--(\ref{365}) all together, we finally obtain that
\begin{equation}\label{366}
\frac{({k_0}\pi)^2}{2L}\mathcal{K}_2=\frac{L((\frac{{k_0}\pi}{L})^2+\lambda)^3\big(\frac{(14-2a_1a_2)(\frac{{k_0}\pi}{L})^2}{(1-a_1a_2)}-\lambda\big)}{48\bar u^2}d_1+O(1).
\end{equation}
Note that $\frac{L((\frac{{k_0}\pi}{L})^2+\lambda)^3}{48{\bar{u}}^2}d_1>0$ in light of (\ref{14}).  Therefore, for $d_1$ and $d_2$ being sufficiently large, we have that
\begin{equation}\label{367}
\text{sgn}\mathcal{K}_2=\text{sgn} \Big(\frac{(14-2a_1a_2)(\frac{{k_0}\pi}{L})^2}{(1-a_1a_2)}-\lambda\Big).
\end{equation}
The rest part of Theorem \ref{theorem32} follows from (\ref{367}) and Proposition \ref{proposition2}.
\end{proof}

According to Proposition \ref{proposition2} and Remark \ref{remark3}, if $L$ is sufficiently small, $k_0=1$ and the only stable bifurcating solution must be spatially monotone.  Theorem \ref{theorem32} further implies $\mathcal K_2>0$ if $L$ is small and $d_1, d_2$ being comparably large, hence the monotone solution must be stable, at least when $\chi$ is slightly larger than $\chi_0$.  On the other hand, when the interval length $L$ is large or the chemical decay rate $\lambda$ is large, the small amplitude bifurcating solution $(u_k(s,x),v_k(s,x),w_k(s,x))$ may be unstable.  Therefore, it is natural to expect the formation of stable steady states of (\ref{11}) with large amplitude, such as boundary layer, interior spikes, etc., or the formation of time periodic spatial patterns.  Rigourous analysis for these spiky solutions requires nontrivial mathematical tools and it is out of the scope of this paper.  Numerical simulations are presented to illustrate the emergence of stable spiky patterns.

\section{Numerical simulations}\label{section4}
In this section, we perform extensive numerical studies to demonstrate the spatial--temporal behaviors of the IBVP (\ref{11}).  To study the effect of chemotaxis on the dynamics of (\ref{11}), we fix $a_1=a_2=0.5$ and $\mu_1=\mu_2=1$, and choose the initial data to be small perturbations of $(\bar u,\bar v,\bar w)$ in all simulations.  Then we choose different sets of parameters to study the regime under which spatially inhomogeneous positive solutions can develop into stable or time--periodic patterns.  The numerical simulations support our theoretical findings on the existence and stability of the bifurcating solutions in Theorems \ref{theorem31} and \ref{theorem32}.  We choose the mesh size of space and time to be $\Delta x = 0.01$ and $\Delta t = 0.01$ respectively in all figures.

\begin{figure}[h!]
        \centering
\includegraphics[width=\textwidth,height=1.5in]{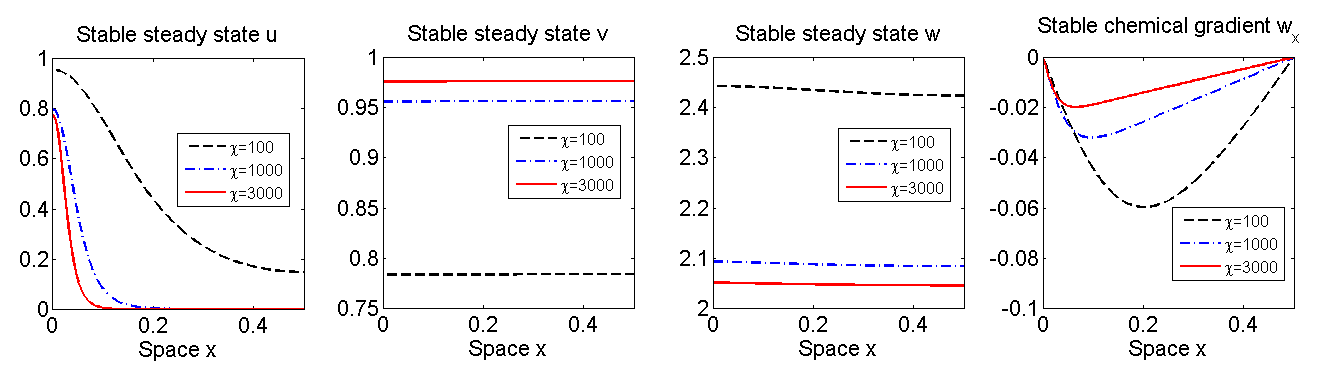}
  \caption{Stable patterns of $u$, $v$, $w$ and $w_x$ of (\ref{11}) for different chemotaxis coefficients $\chi$.  The parameters are $d_1=1$, $d_2=0.1$, $a_1=a_2=0.5$, $\mu_1=\mu_2=1$, and $\lambda=\xi=0.5$.  The initial data are small perturbations from the homogeneous equilibrium $(u_0,v_0,w_0)=(\bar u,\bar v,\bar w)+(0.01,0.01,0.01)\cos (2.4\pi x)$.  These graphes suggest that small interval supports stable monotone steady states and large chemotaxis drives the formation of stable spiky solution.  Moreover, the positive steady state $u$ concentrates and shifts to the boundary as $\chi$ increases.}\label{fig2}
\end{figure}

Table \ref{table1} lists the values of $\chi_k$ and $\hat \chi_k$ for $k=1,2,...,7$ when $L=0.5$.  We see that $\min_{k\in\mathbb N^+}\{\chi_k,\hat \chi_k\}=\chi_1=61.0$, therefore the stable steady state of (\ref{11}) must be spatially monotone according to our stability theorem, at least when $\chi$ is slightly larger than 61.  This is illustrated through steady state $u(x)$ given in Figure \ref{fig2} where $\chi$ is taken to be 100, 1000 and 3000.  We see that the stable steady state $u$ has a single boundary layer at $x=0$ which shifts to a boundary spike as the chemotaxis coefficient $\chi$ increases.  Steady state $v$ is almost constant since its chemotaxis coefficient $\xi$ and the chemical gradient $w_x$ is very small.  The simulation indicates that chemotaxis is a dominating mechanism for the formation of stable patterns in (\ref{11}), since $u$ and $v$ have the same kinetics but quite different patterns.
\begin{table}[h!]
\centering
\begin{tabular}{|c|c|c|c|c|c|c|c|} \hline
$k$& 1 & 2 & 3 & 4 & 5& 6& 7\\ \hline
$\chi_{k}$&\textbf{61.0}& 238.6 & 534.7& 949.2& 1482.2&  2133.6 & 2903.4 \\ \hline
$\hat \chi_{k}$&75.2& 290.2 & 648.4&  1150.0& 1794.9&  2583.1 &  3514.6 \\
\hline
\end{tabular}
\caption{List of values of $\chi_k$ in (\ref{39}) and $\hat \chi_k$ in (\ref{35}) for $L=0.5$.  The parameters are $d_1=1$, $d_2=0.1$, $a_1=a_2=0.5$, $\mu_1=\mu_2=1$, and $\lambda=\xi=0.5$. We see that $\min_{k\in\mathbb N^+}\{\chi_k,\hat\chi_k\}=\chi_1$.  Therefore $(\bar u,\bar v,\bar w)$ lost its stability to the stable mode $\cos \frac{\pi x}{L}$ which is monotone in space.  This is numerically verified in Figure \ref{fig2}.}
\label{table1}
\end{table}

\begin{figure}[h!]
\centering
\includegraphics[width=\textwidth,height=2.5in]{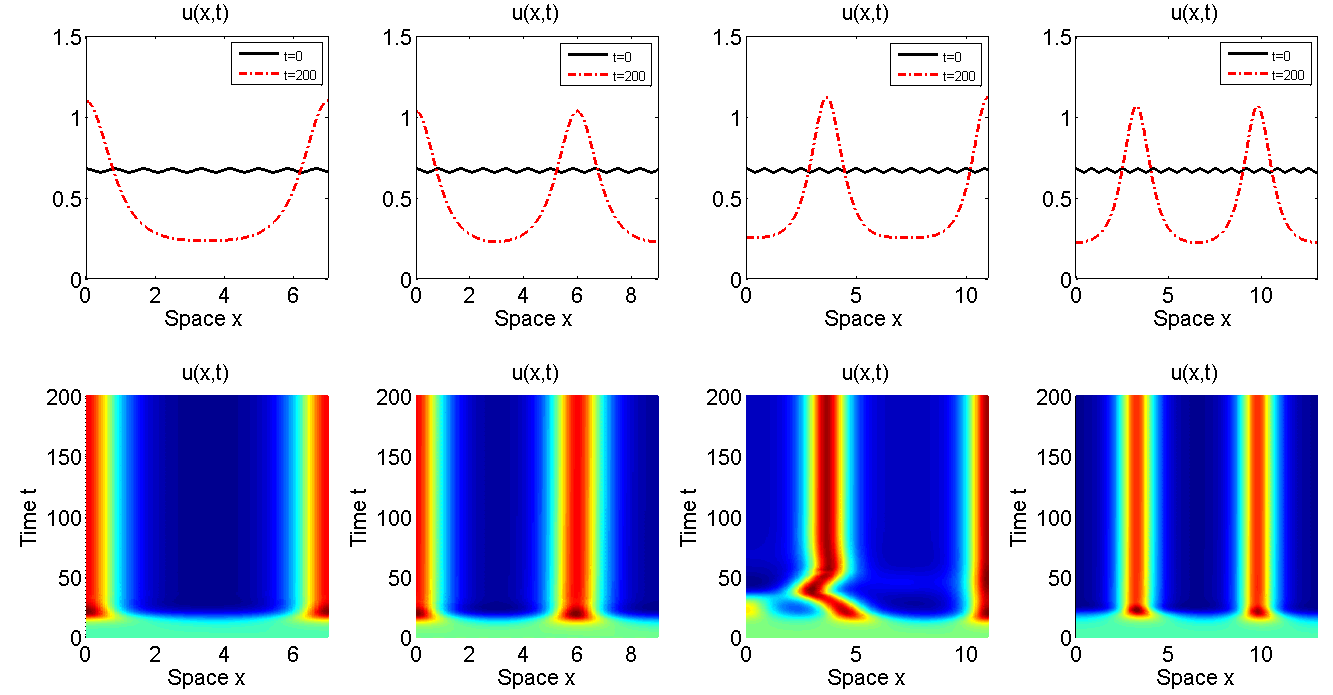}
  \caption{Formation of stable patterns over intervals with different lengthes $L=7$, $9$, $11$ and $13$.  The system parameters and initial data in all graphes are taken to be the same as in Figure \ref{fig2} except that $\chi=6$, which is slightly larger than $\chi_{k_0}$ given in Table \ref{table2}.  These graphes support our stability analysis of the bifurcating solutions.}\label{fig3}
\end{figure}

In Table \ref{table2}, we list the values of minimum steady state bifurcation parameter $\chi_0$ and the corresponding wavemode number $k_0$ for the interval length $L=3,5,...,21$.  In Figure \ref{fig3}, we plot the spatial--temporal behaviors of $(u(x,t),v(x,t),w(x,t))$ to demonstrate the stable wavemode selection mechanism for $L=7,9,11,13$, where $k_0=2,3,3,4$.  According to our stability results in Theorem \ref{theorem32}, $(\bar u,\bar v,\bar w)$ loses its stability to the stable wavemodes $\cos \frac{2\pi x}{7}$, $\cos \frac{3\pi x}{9}$, $\cos \frac{3\pi x}{11}$ and $\cos \frac{4\pi x}{13}$ respectively.  Figure \ref{fig3} verifies such stable wavemode selection mechanism for (\ref{11}) where we take $\chi=6$ to be slightly larger than $\chi_0$.
\begin{table}[h!]
\centering
\begin{tabular}{|c|c|c|c|c|c|}
  \hline
Interval length $L$& 3 & 5 & 7& 9& 11\\ \hline
    $k_0$&1&1&2&3&3\\ \hline
    $\chi_0=\chi_{k_0}$&4.5868& 4.8455& 4.4330& 4.5868& 4.4260\\ \hline
Interval length $L$& 13 & 15 & 17 & 19& 21\\ \hline
    $k_0$&4&4&5&5&6\\ \hline
    $\chi_0=\chi_{k_0}$&4.4815& 4.4290& 4.4465&  4.4327& 4.4330\\
  \hline
\end{tabular}
\caption{List of stable wavemode numbers and the corresponding bifurcation values for different interval lengthes.  System parameters are chosen to be the same as those in Table \ref{table1}.  We see that larger intervals support higher wavemodes.}
\label{table2}
\end{table}

\begin{figure}[h!]
        \centering
\includegraphics[width=\textwidth,height=2.5in]{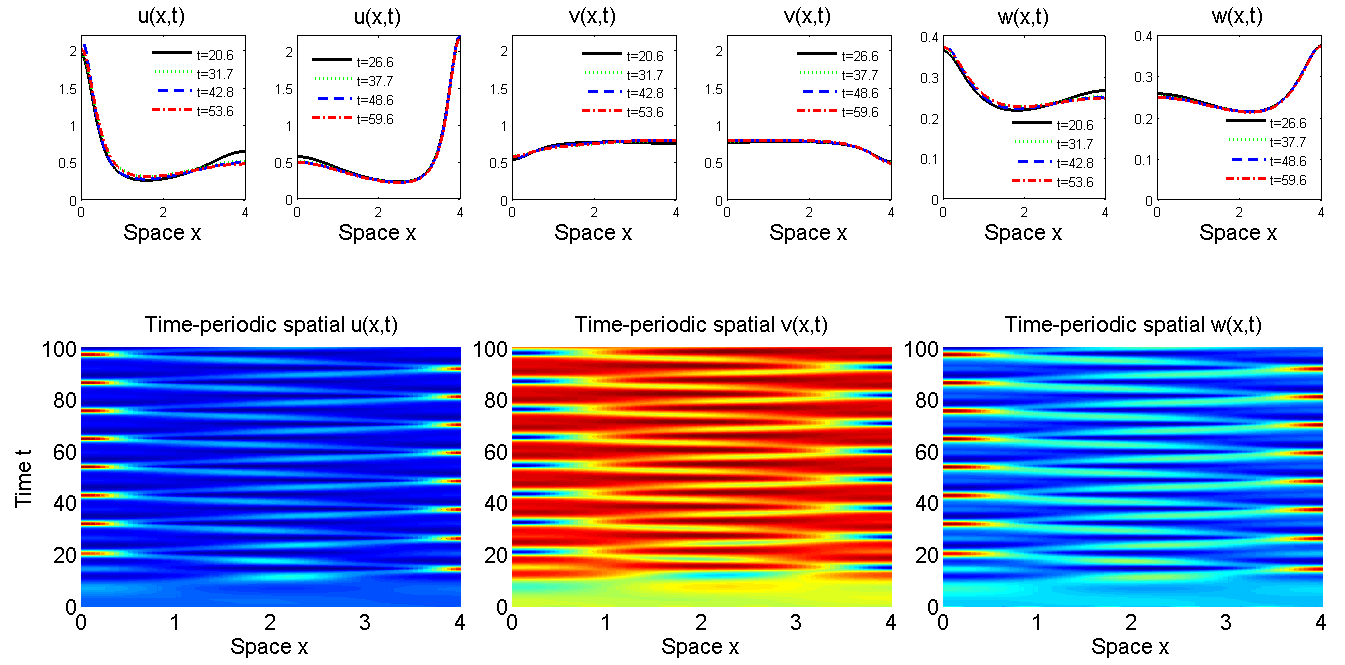}
  \caption{The formation and evolution of time--periodic spatial patterns over $\Omega=(0,L)$.  The parameters are $d_1=5$, $d_2=0.1$, $a_1=a_2=0.5$, $\mu_1=\mu_2=1$, and $\lambda=5$, $\xi=0.1$.  $\chi=80$ is slightly larger than the least bifurcation value 66.4 for $L=4$ as shown in Table \ref{table3}.  The initial data are small perturbations from the homogeneous equilibrium $(u_0,v_0,w_0)=(\bar u,\bar v,\bar w)+(0.01,0.01,0.01)\cos (2.4\pi x)$.  Our simulations demonstrate that $u$, $v$ and $w$ are all time--periodic and this suggests that Hopf bifurcation occurs at $\chi_0=\hat \chi_4$.}\label{fig4}
\end{figure}

\begin{figure}[h!]
        \centering
\includegraphics[width=\textwidth,height=2.5in]{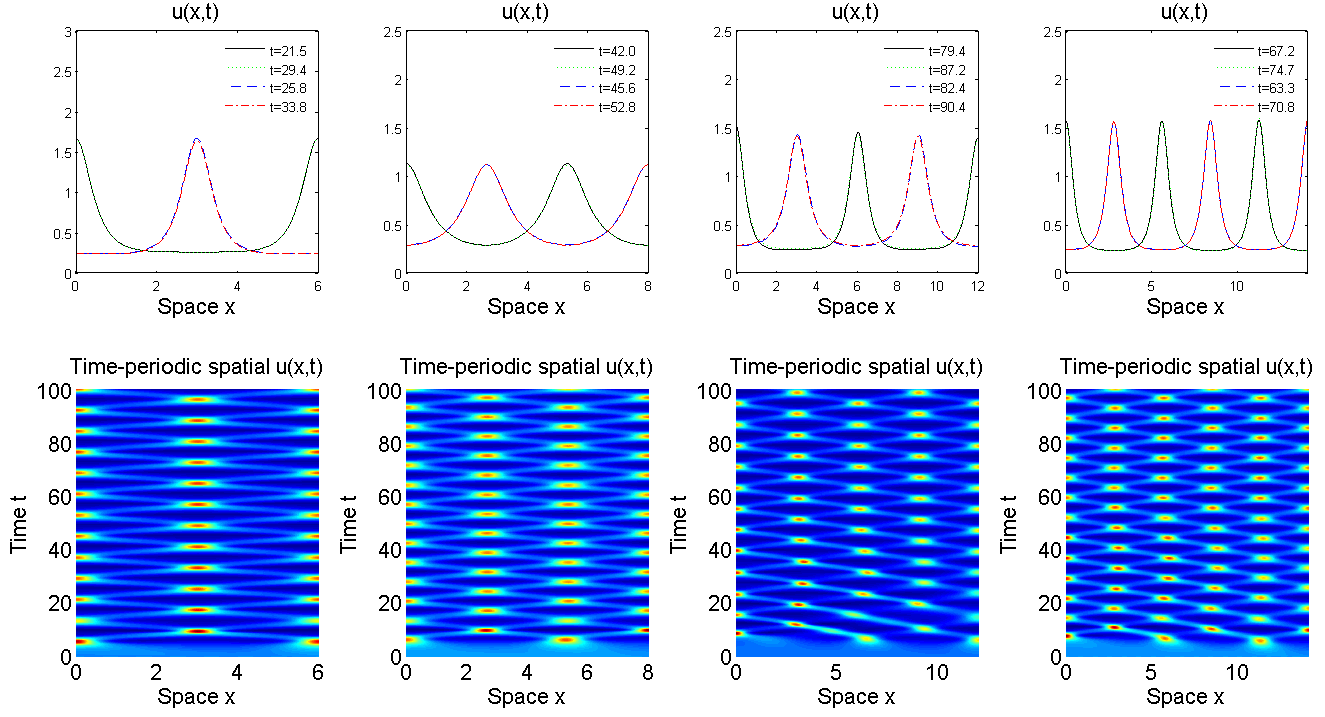}
  \caption{Spatial--temporal behaviors of $u(x,t)$ of (\ref{11}) over intervals with different lengthes.  Parameters and initial data are chosen to be the exact same as in Figure \ref{fig4}.  Our simulations indicate that the periodic pattern has a stable wavemode $\cos \frac{k_1 \pi x}{L}$ where $k_1$ is given in Table \ref{table3} for different intervals.}\label{fig5}
\end{figure}

Figure \ref{fig4} and Figure \ref{fig5} illustrate the emergence and evolution of time--periodic spatial patterns of (\ref{11}).  Taking $d_1=5$, $d_2=0.1$, $a_1=a_2=0.5$, $\mu_1=\mu_2=1$, and $\lambda=5$, $\xi=0.1$, we find that $\chi_0$ is always attained by $\hat \chi_{k_1}<\min_{k\in \mathbb N^+} \chi_k$ as listed in Table \ref{table3}.  According to Case (ii) in Proposition \ref{proposition2} or Theorem \ref{theorem32}, all the steady state bifurcating solutions are unstable around the bifurcating points $(\bar u,\bar w,\bar w, \chi_k)$, hence we surmise that Hopf bifurcation occurs at $\chi=\chi_{k_1}$.  In Figure \ref{fig4}, we plot the spatial--temporal behaviors of solutions $(u(x,t),v(x,t),w(x,t))$ to (\ref{11}) which are time--periodic spatially inhomogeneous.  The period is approximately $T=11$.  In Figure \ref{fig5}, we demonstrate the effect of interval length on the wavemode of the periodic patterns.  According to Table \ref{table3}, $k_1=2,3,4,5$ when $L=6,8,12,14$, and we surmise that the periodic pattern has a stable wavemode $\cos \frac{k_1 \pi x}{L}$ to which $(\bar u,\bar v,\bar w)$ loses its stability at $\chi_0=\hat \chi_{k_1}$.  Moreover, we can find that the periods are approximately $T=7.85,7.20,7.88$ and $7.47$ in each plot.  Rigourous analysis of the time--periodic spatial patterns is beyond the scope of this paper.  It is also interesting and important to study the effect of domain size or system parameters on the period.
\begin{table}[h!]
\centering
\begin{tabular}{|c|c|c|c|c|c|c|c|}
  \hline
Interval length $L$& 1 & 2 & 3& 4& 5& 6& 7\\ \hline
$k_1$&1 & 1 & 1& 1& 2& 2& 2\\ \hline
$\chi_0=\hat \chi_{k_1}$&129.0 &  69.7 & 63.2& 66.5& 64.5& 63.2& 64.2\\ \hline
Interval length $L$& 8 & 9 & 10 & 11& 12& 13 & 14\\ \hline
    $k_1$&3 & 3 & 3& 4& 4& 4& 5\\ \hline
    $\chi_0=\hat \chi_{k_1}$&63.8 & 63.2 & 63.7& 63.5& 63.2& 63.5& 63.4\\
  \hline
\end{tabular}
\caption{List of values of $\chi_0$ for different interval lengthes.  The system parameters are $d_1=5$, $d_2=0.1$, $a_1=a_2=0.5$, $\mu_1=\mu_2=1$, and $\lambda=5$, $\xi=0.1$.  In this case, $\chi_{0}$ is achieved at $\hat \chi_{k_1}$ and our numerical simulations in Figure \ref{fig4} and Figure \ref{fig5} suggest that Hopf bifurcation occurs at $\chi=\hat \chi_{k_1}$ and the wavemode $\cos \frac{k_1 \pi x}{L}$ is stable.}
\label{table3}
\end{table}
\begin{figure}[h!]
        \centering
\includegraphics[width=\textwidth,height=2.5in]{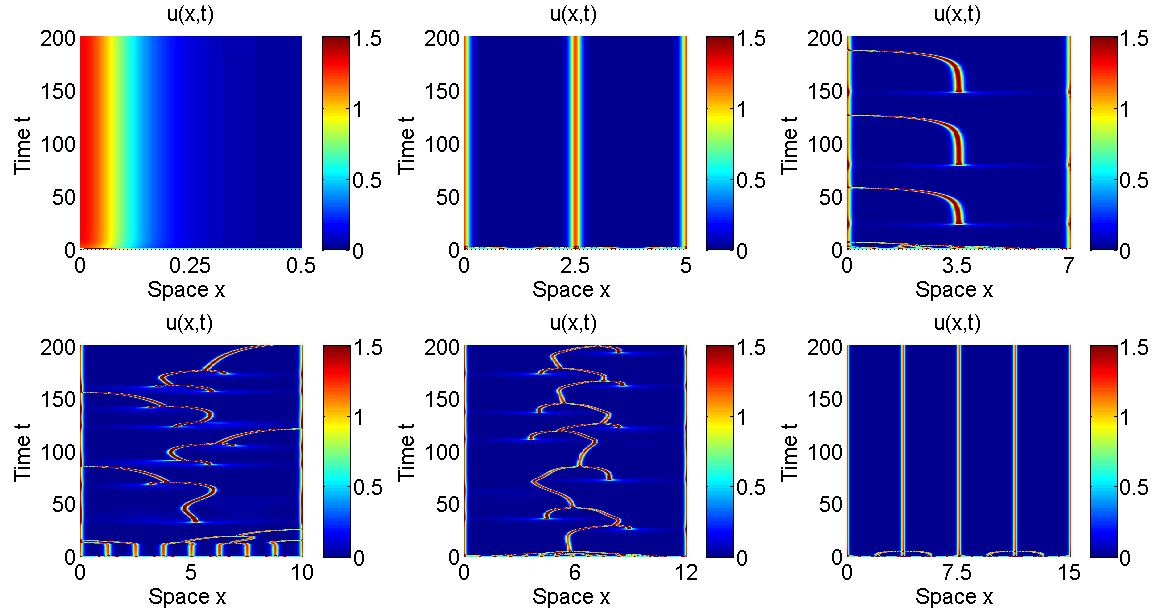}
  \caption{Formation of stable spikes and coarsening patterns of $u(x,t)$.  We choose $d_1=0.2$, $d_2=0.3$, $\lambda=0.5$ and $\chi=20$, $\xi=50$, while the rest parameter are the same as in Figure \ref{fig4}.  The initial data are $(u_0,v_0,w_0)=(\bar u,\bar v,\bar w)+(0.01,0.01,0.01)\cos (2.4\pi x)$.  We see that the positive solution develops into a stable single or multi--spike when $L=0.5$, $5$ and $15$, while coarsening occurs with the merging and emerging of spikes when $L=7,10$ and $12$.}\label{fig6}
\end{figure}
Finally, we present in Figure \ref{fig6} a set of simulations of the dynamical behaviors of spatially inhomogeneous patterns of (\ref{11}) as domain size $L$ changes, while $\chi$ is taken to be far away from the bifurcation values.  Similar as in our wavemode selection mechanism of the bifurcating solutions, we find that the large domain size $L$ tends to increase the number of spikes or cell aggregates.  It is also observed that multi--spiky solutions undergo a coarsening process in which two interior spikes merge into a single spike and this new spike develops and merges with another spike.  Recently the authors constructed time--monotone Lyapunov functional to (\ref{11}) in \cite{WYZ} when $\mu_1=\mu_2=0$, therefore the cellular growth is necessary for the formation of time--periodic or coarsening patterns to (\ref{11}).  Moreover the instability of the interior spikes is needed to investigate such spike merging phenomenon, however in contrast to the classical Keller--Segel model where the interior spike is known to be unstable, some interior spikes are stabilized for some proper length $L$, for example when $L=15$.  However, rigorous analysis is required if one wants to find the critical threshold responsible for spike merging and stabilization.

\section{Conclusions and Discussions}\label{section5}

In this paper, we study the $3\times 3$ chemotaxis system (\ref{11}) of two competing species $u$ and $v$ with Lotka--Volterra type kinetics.  It is assumed that both species move chemotactically along the concentration gradient of the same chemical $w$ which is also produced by the species.  This system is the one--dimensional counterpart of (\ref{12}) with $\tau=1$ which was studied by Tello and Winkler in \cite{TW}.  When $\tau=0$, they showed that the homogeneous equilibrium $(\bar u,\bar v,\bar w)$ is the global attractor of (\ref{12}) if the chemo--attraction rates $\chi$ and $\xi$ are small compared to the intrinsic growth rates $\mu_1$ and $\mu_2$ as in (\ref{15}).

We obtain the global existence and uniform boundedness of classical positive solutions $(u(x,t),v(x,t),w(x,t))$ of the one--dimensional system (\ref{11}).  For the stationary system (\ref{31}), we show that the homogeneous equilibrium $(\bar u,\bar v,\bar w)$ becomes unstable for $\chi>\chi_0=\min_{k\in \mathbb N^+}\{\tilde \chi_k, \hat \chi_k\}$ given by (\ref{33}) in the sense of Turing's instability driven by chemotaxis or advection.  Then we apply Crandall--Rabinowitz bifurcation theories to establish the existence and stability of its nonconstant positive steady states that bifurcate from $(\bar u,\bar v,\bar w, \chi_k)$ with $\chi_k=\tilde \chi_k$.  We show that the bifurcation branches are of pitch--fork type and perform detailed calculations to obtain the formula of constant $\mathcal K_2$ that determines the turning direction of each branch around the bifurcation point.  Our main results provide a selection mechanism of stable wavemode for system (\ref{11}).  If $\chi_0=\chi_{k_0}<\min_{k \in \mathbb N^+} \hat \chi_k$, then $(\bar u,\bar v,\bar w)$ can only lose its stability to the stable steady state bifurcating solution with wavemode $\cos \frac{k_0 \pi x}{L}$, at least when $\chi$ is around but larger than $\chi_0$.  If $\chi_0=\hat \chi_{k_1}<\min_{k \in \mathbb N^+}\chi_k$, then all bifurcating solutions are unstable independent of the turning direction of the local branches.  Our stability results suggest that when the interval length $L$ is sufficiently small, the stable steady state of (\ref{11}) must be spatially monotone.

Numerical simulations of (\ref{11}) are performed to support our theoretical findings about the stable wavemode selection mechanism.  Our numerics also illustrate the dynamical behaviors of solutions to (\ref{11}) that exhibit complex spatial--temporal patterns such as merging and emerging of spikes, time--periodic spatial patterns with interesting structures, etc.  It turns out the chemotaxis strength and interval length are very important in the formation and evolution of these interesting patterns of (\ref{11}).

There are a few questions that we want to propose for further studies.  First of all, it is interesting and also important to study the global existence of (\ref{12}) over high--dimensional spaces, when the assumption (\ref{15}) is relaxed, in particular over $\Omega \subset \mathbb{R}^2$.  The literature suggests that blow--ups can be inhibited by the degradation in the cellular kinetics, however, whether or not this is sufficient to prevent blow--ups in finite or infinite time over high--dimensional spaces when chemo--attraction rates $\chi$ and $\xi$ are large remains open so far.  Both $\chi$ and $\xi$ are required to be positive in the arguments of \cite{TW}, and in light of our stability analysis, we surmise that the solution of (\ref{12}) is always global if $\chi<0$ and $\xi<0$, since chemo--repulsion has smoothing effect like diffusions.  To prove this, one needs an approach totally different from that in \cite{TW}.

Our stability analysis of the homogeneous solution $(\bar u,\bar v,\bar w)$ and numerical simulations of the time--periodic spatial patterns in Figure \ref{fig4} and Figure \ref{fig5} suggest that the homogeneous equilibrium loses its stability through Hopf bifurcation if $\chi_0=\hat \chi_{k_1}<\min_{k\in\mathbb N^+}\tilde \chi_k$.  Moreover, the stable Hopf bifurcating solution has the wavemode number $k_1$ which is increasing or at least non-decreasing in $L$.  However, rigorous analysis is needed for the existence and stability of these periodic solutions.  It is also interesting to study how the period depends on the chemotaxis rate and domain size.

Theoretical analysis of the existence and stability boundary and interior spikes in Figure \ref{fig2} and Figure \ref{fig6} is another challenging problem worth giving attention to, even in one--dimensional domain.  It is well known that large chemotaxis rate supports the formation of spiky solutions in $2\times2$ system.  Whether or not this is true for $3\times 3$ system is an important but very challenging problem that one can pursue in the future, in particular when Hopf bifurcation occurs.  Moreover, the stability of spiky solutions to (\ref{12}) also deserves future exploring.  We also want to mention that $d_1$ and $d_2$ are required large in evaluating the sign of $\mathcal K_2$, however, we need some extremely complicated and difficult calculations in order to remove this constraint.

\end{document}